\input amstex
\documentstyle{amsppt}
\magnification = 1200
\loadbold
\NoRunningHeads
\pageno=1
\expandafter\let\csname logo\string@\endcsname=\empty
\pagewidth{5.5in}
\pageheight{7.25in}
\NoBlackBoxes

\define\a{\alpha}
\predefine\barunder{\b}
\redefine\b{\beta}

\predefine\dotunder{\d}
\redefine\d{\delta}

\define\f{\frac}
\define\g{\gamma}

\define\lb\{{\left\{}

\define\lra{\longrightarrow}

\define\om{\omega}

\define\oper{\operatorname}

\define\p{\partial}
\define\rb\}{\right\}}
\define\s{\sigma}

\define\th{\theta}

\define\wt{\widetilde}

\define\({\left(}
\define\){\right)}
\define\[{\left[}
\define\]{\right]}
\define\<{\left<}
\define\>{\right>}
\def\slantline#1#2#3#4#5{\hbox to 0pt{
}}


\def\SC{{\Cal C}}
\def\SD{{\Cal D}}
\def\SE{{\Cal E}}
\def\SF{{\Cal F}}

\def\SI{{\Cal I}}
\def\SJ{{\Cal J}}
\def\SK{{\Cal K}}

\def\SO{{\Cal O}}

\def\SR{{\Cal R}}
\def\SS{{\Cal S}}

\def\SU{{\Cal U}}

\def\SW{{\Cal W}}


\def\BD{{\bold D}}
\def\BE{{\bold E}}
\def\BF{{\bold F}}

\def\BL{{\bold L}}


\def\bbp{{\Bbb P}}

\def\bbr{{\Bbb R}}

\def\bbz{{\Bbb Z}}

\define\od{oriented }

\define\loc{L^1_{\text{loc}}}

\define\Hom{\oper{Hom}}

\define\Div{\oper{Div}}

\define\bsnu{\boldsymbol\nu}

\define\UP{{U^\bot}}

\define\ur{\underline \bbr}

\define\Ua{{\Cal U}_\alpha}

\define\OV{{\Cal O}_{V}}

\define\orX{{\Cal O}_{X}}

\define\pxtx{$p: {\wt X}\to X\,\,$}
\define\Vt{{\wt V}}
\define\Xt{{\wt X}}

\define\orV{{\Cal O}_{V}}

\define\orTX{{\Cal O}_{X}}

\define\orF{{\Cal O}{_F}}

\define\stF{V_{n-q+1} (F)}
\define\pist{\pi_{q-1}(\stF )}

\define\LOC{{\oper{loc}}}
\define\nut{\wt{\nu }}

\topmatter
\title
Stiefel--Whitney  Currents 
\endtitle
\author Reese Harvey and  John Zweck \endauthor

\abstract\nofrills{\usualspace }

A canonically defined {\it mod 2 linear dependency current} is
 associated to each collection $\nu$ of sections, $\nu_1,...,\nu_m$,
 of a real rank $n$ vector bundle. This current is supported on the
 linear dependency set of $\nu$. It is defined whenever the 
 collection $\nu$ satisfies a weak measure theoretic condition called 
\lq\lq atomicity''. Essentially any reasonable collection of sections
 satisfies this condition, vastly extending the usual general
 position hypothesis. 
This current is a mod 2 $d$--closed locally integrally flat current
 of degree $q = n-m+1$ and hence determines a  $\bbz_2$--cohomology 
class. This class is shown to be well defined independent of the
 collection of sections. Moreover it is the $q$th Stiefel--Whitney
 class of the vector bundle. 

More is true if $q$ is odd or $q=n$. In this case a
{\it linear dependency current}
which is twisted by the orientation of the bundle
  can  be associated to the collection
 $\nu$. The mod 2 reduction of this current is the mod 2 linear
 dependency current. The cohomology class of the linear dependency
 current is 2-torsion and  is the $q$th twisted integral Stiefel--Whitney
 class of the bundle.

 In addition, higher dependency and
general degeneracy currents of bundle maps are studied, together with
 applications to singularities of projections and maps.

These results rely on a theorem of Federer which states that the complex
 of integrally  flat currents mod $p$ computes cohomology mod $p$. 
An alternate  approach to Federer's theorem
 is offered in an appendix. This approach is simpler and is via sheaf theory.

\endabstract
\address  Rice University, Houston, Texas.\endaddress
\address University of Nevada, Reno, Nevada. \endaddress
\endtopmatter
 
\document

\subheading{Contents}
\newline\newline
1. Introduction.\newline
2. Divisors and atomicity.\newline
3. Linear dependency currents.\newline
4. Stiefel--Whitney currents.\newline
5. Obstruction currents.\newline
6. Higher dependency currents.\newline
7. Applications\newline
A. Computing cohomology with currents.
\newline

\subheading{1. Introduction}
\newline
\newline

It is well known [BC,Pon,S,W] that the linear dependency set of a
 collection of sections of a vector bundle is related to the characteristic 
classes of the bundle. 
In particular the zero set of a regular section 
defines a cohomology class which is the Chern-Euler class  of the bundle.
 In [HL3]  canonically  defined current 
representatives of the Chern classes of a complex vector bundle were 
associated to collections of smooth sections of the bundle.
 These currents are called 
 linear dependency currents since they are supported on the linear
 dependency set of the collection of sections. 

The main aim of this paper is to study the linear dependency currents of a
 collection of sections of a {\it real} vector bundle.
 These are either mod 2 or
 bundle-twisted currents which represent either the mod 2 or twisted-integer
Stiefel-Whitney classes of the bundle. Since they are either mod 2 or 2-torsion
these currents were overlooked in [HL3].

The linear dependency current associated with an ordered collection $\nu$ of
sections of a real vector bundle is defined in Section 3 paralleling 
 a standard 
construction in enumerative geometry (see for example [Por], [Fu]). 
In general  the most one can say is that this linear dependency current, $\BL\BD
^{\oper{mod 2}} (\nu)$,
 is a mod 2 current and that it  determines a $\bbz_2$-cohomology class which
is well defined independent of the particular collection of sections 
of the bundle, (Theorem 3.15). 

However if the degree of  $\BL\BD^{\oper{mod 2}} (\nu)$ is odd or equal to the
rank of the bundle it is also possible to define a (bundle-twisted) linear
dependency current, $\BL\BD (\nu)$, which encodes certain (twisted) integer
multiplicities of dependency among the sections, (Proposition 3.17). 
The mod 2 reduction of this current $\BL\BD (\nu)$ is the mod 2
 current  $\BL\BD^{\oper{mod 2}} (\nu)$. The current $\BL\BD (\nu)$ determines
a (twisted) integer cohomology class well defined independent of the choice
of collection $\nu$, (Theorem 3.7). If the degree of $\BL\BD (\nu)$ is less
than the rank of the bundle (which occurs when the collection consists
of more than one section) this cohomology class is a torsion class of order
2, (Corollary 3.9).

A major advantage of the approach taken here is that the linear dependency
current is defined whenever the collection $\nu$ satisfies a weak measure
theoretic condition called \lq\lq atomicity\rq\rq, which is vastly more
 general than the usual general position hypothesis. For example, a real 
analytic
collection of $m$ sections of a rank $n$ bundle is atomic provided that, 
for all $j\in \{ 0,1,...,m-1 \}$, the codimension of the set of points 
over which exactly $j$ of the sections are linear independent is at least
the expected codimension $n-j$, (see [HL3,2.14]). 

Another important property of the (mod 2) linear dependency current is that
 it is a (mod 2) locally integrally flat current. Recall that the integrally
 flat currents are those of the form $R+dS$, where $R$ and $S$ are rectifiable.
Federer [F] proved that the complex of locally integrally flat currents 
(or such currents mod $p$) can be used to compute integer (or mod $p$)
cohomology. In the Appendix we offer an alternate approach to the theory of 
(mod $p$) integrally flat currents and their cohomological properties.
This simple approach is via sheaf theory and is quite distinct from the
form of the theory given in the literature.

The theory of dependency currents relies heavily on the theory of 
zero divisor currents which was originally
 developed  in [HS] for
 \lq\lq atomic \rq\rq sections of an oriented vector bundle 
over an oriented manifold.
The notion of an atomic section provides a generalization of the notion of a 
section being transverse to zero, one which is both useful and vastly more general. The zero divisor is a $d$--closed locally integrally flat 
 current which determines a  unique
integer cohomology class, the Euler class. In this paper it is crucial that the
notion of a zero current be understood in the non--orientable case. This is done in Section 2 where the zero divisor is defined as a bundle  twisted current. This current determines a cohomology class, (Theorem 2.5), which is the twisted
Euler class, $\wt e \in H^n(X,\wt\bbz )$, of the vector bundle. The reduction mod 2 of the zero divisor eliminates the twisting, yielding a mod 2 current
 which represents the top Stiefel--Whitney class, $w_n \in H^n(X,\bbz_2)$.

In Section 4 we identify the
 $\bbz_2$--cohomology class of the
degree $q$  mod 2 current $\BL\BD^{\oper{mod} 2} (\nu)$ as $w_q$, the $q$th Stiefel--Whitney class of the bundle,  (Theorem 4.1). Moreover,  if $q$ is odd, the $\wt\bbz$--class of the twisted current $\BL\BD (\nu)$ is identified as the $q$th twisted integral Stiefel--Whitney class $\wt W_q\in H^q(X,\wt\bbz)$, (Theorem 4.10). This result
 is a corollary of the fact that the Bockstein of the mod 2 dependency current $\BL\BD^{\oper{mod} 2} (\nu)$
of degree $q-1$ is the degree $q$ twisted dependency current associated with a subcollection of the collection $\nu$, ( c.f. [W], [St]).

The Stiefel--Whitney classes were originally defined ([S], [W], [St]) as the 
primary obstruction to the existence of certain collections of linearly independent sections of a bundle $F\to X$. In Section 5 we examine the relationship between linear dependency currents and obstruction cocycles. Given a triangulation of $X$ it is possible to choose a particular atomic collection of sections of $F$ so that the Steenrod obstruction cocycle of the collection is defined. The Poincar\'e dual of such a cocycle is a cycle which defines a current on $X$ by integration. We then show that this obstruction current is equal to the linear dependency current of the particular collection of sections.
 Among other things, this provides an alternate proof of the results of Section 4. 

In Section 6 higher dependency currents and general degeneracy currents of
vector bundle maps are discussed, further expanding  the results of [HL3].
Some of the degeneracy currents studied in Section 6 were not included in [HL3] since they are either not defined as twisted currents or their real cohomology
class is zero. In these cases we can define mod 2 and/or twisted  degeneracy
 currents. The integer cohomology classes of the twisted degereracy currents
 were
first studied by Ronga [R] who proved that they are uniquely determined by 
their torsion free part and mod 2 reduction. We expand upon Ronga's 
result by explicitly identifying the integer cohomology classes 
of the higher dependency currents as  certain polynomials 
in the integer Pontrjagin and Stiefel--Whitney classes, (Theorem 6.15).

In Section 7 applications of the theory to singularities of projections and
maps are given. In particular we recover the well know fact that the Steifel--Whitney classes of the tangent bundle $TX$ and normal bundle $NX$ of a submanifold $X\subset \bbr^N$ can be defined in terms of  singularities of projections. The original version of this result is  due to 
 [Pon], [T] (see also [BMc]). Note however that they only consider generic projections whose critical sets are non--degenerate, with multiplicity $\pm 1$.
The atomic theory enables us to consider degenerate critical sets of arbitrary integer multiplicity,  (Proposition 3.17).
 In particular, if $X$ is a real analytic submanifold, the tangent and normal Stiefel--Whitney classes can be defined in terms of the singularities of any projection whose  degeneracy subvarieties have at least  the expected codimension.
Integer and mod 2 cohomological obstructions to the existence 
of smooth immersions and surjections between manifolds are also
 given, c.f. [R].

Two further applications are worth noting. Following [HL3] we can define mod 2 and twisted integer degeneracy currents associated with higher self intersections of plane fields and invariants of pairs of foliations. Mod 2 and integer umbilic currents of hypersurfaces can also be studied using these ideas. Details of these two applications are left to the reader. 

Secondary (Cheeger-Chern-Simons) Stiefel-Whitney classes will be 
introduced in a later paper. Canonical $\loc$ representatives of these
classes will be associated to each atomic collection of sections of a bundle
with riemannian connection. In the case of a single section $\a$ the secondary Euler class is represented by the Chern-Euler potential $\s (\a )$. 
As is discussed in [HL1] this potential satisfies the important equation
$d\s (\a) \,\,=\,\,\chi\,\,-\,\,\Div (\a)$, where $\chi$ is the Euler form
and $\Div (\a)$ the divisor of the section. If the collection $\nu$ consists of
more than one section then there is a canonical $\loc$ current $T(\nu )$,
satisfying the current equation
$ dT(\nu) \,\,=\,\, -\, \BL\BD (\nu)$, which represents
 the appropriate secondary
Stiefel-Whitney class.
This current equation is related to a formula of Eells, [E].

\subheading{2. Divisors and atomicity}
\newline
\newline

Harvey and Semmes defined the zero divisor current of an atomic section of an {\it oriented} real rank $n$ vector bundle over an {\it \od } manifold. The divisor is a 
codimension $n$ current which is supported on the zero set of the section and which encodes the integer multiplicity of vanishing of the section. Furthermore 
it is a $d$--closed locally integrally flat current whose cohomology class in $H^n(X,\bbz )$ is well defined independent of the choice of section. This class is the Euler class of the bundle.

The aim of this section is to define and study the zero divisor current in the case in which neither the vector bundle nor the base manifold are assumed to
 be orientable. In this case the zero divisor is defined to be a current which is twisted by the orientation bundle of the vector bundle. It is also useful
to define the mod 2 divisor to be the mod 2 reduction of the divisor. Both
of these notions of divisor will be important  in the study of linear
dependency currents in Section 3.

We begin by recalling some definitions. Let $V\to X$ be a real rank $n$ vector
bundle over an $N$-dimensional manifold. No orientation assumptions will
be made on $V$ or $X$. Let $\orX$ and $\orV$ denote the principal
 $\bbz_2$-bundles
of orientations of $TX$ and $V$ over $X$. An $\Cal O${\bf-twisted $k$-form}
is a section of $\Cal O\otimes_{\bbz_2}\Lambda^kT^*X\to X$. (Often the subscript $\bbz_2$ will be dropped when tensoring with $\Cal O$.) 
 A {\bf density}  is a top degree 
smooth $\orX$-twisted form on $X$.
 Note that densities can be integrated over $X$.
A {\bf generalized function} 
is a continuous linear functional on the space of compactly supported smooth 
densities on $X$.
A {\bf current} of degree $p$   is a differential $p$--form on $X$ whose 
coefficients (with respect to each coordinate system) are
 generalized functions.
Equivalently, a degree $p$ current is a continuous linear functional on 
the space of compactly supported $\orX$-twisted $(N-p)$--forms.
Similarly an $\orV${\bf-twisted current} is an $\orV$-twisted  form whose coefficients are generalized functions, i.e. it  acts on $\orV\otimes\orX$-twisted forms. An $\loc$ form is a form whose coefficients are $\loc$ functions. Therefore $\loc$ forms are  currents which are not twisted. On the other hand an oriented compact submanifold of $X$ defines an $\orX$-twisted current by integrating 
(untwisted) forms over it.
 Note that exterior differentiation is a well defined operation on
(twisted) currents.

 On a  contractable 
 open subset $U$ of $X$ each $\OV$--twisted current $\wt T$ can be written in the form $\wt T = [e]\otimes T$ where $[e]\in\SO_V$ is the orientation class of a local frame $e$ for $V$ over $U$ and where $T$ is a 
current on $U$. If $V$ is orientable each choice of orientation defines an isomorphism
between $\orV$-twisted currents  and currents. These two isomorphisms differ by a minus sign. Note that the definition of a current on a non-orientable manifold
agrees with that given in [HL4] but disagrees with that in [Z1,2].

\subheading{A. Divisors in the nonorientable case}
\newline\newline

In this subsection we define and study the divisor of a section of $V \to X$.
 The divisor is defined to be an $\orV$-twisted current. Note that, if $V$ and $X$ are 
oriented, the definition of divisor given below agrees with that of [HS].

 The \bf solid angle kernel, \rm $\theta$, 
is the $\loc$ form on $\bbr^n$ obtained by pulling back the normalized volume form on the unit sphere to $\bbr^n\sim\{ 0\} $ by the radial projection map. The  current equation $d\, \th  = [0]$ on $\bbr^n$, where $[0]$ denotes the point mass at the origin, motivates the definition of divisor.

\subheading{Definition 2.1}
Let $X$ be a smooth manifold and let $y=(y_1,...,y_n)$ denote coordinates on 
$\bbr^n$. 
In the case $n > 1$ a smooth function  $u: X\to \bbr^n$ is called {\bf atomic} if, for each form $\f{dy^I}{|y|^p}$ on $\bbr^n$ with $p = |I| \leq n-1$, the pullback $u^*(\f{dy^I}{|y|^p})$ to $X$ has an $\loc (X)$ extension across the zero set $Z$ of $u$. Also assume that $u$ does not vanish identically in any connected component of its domain $X$. 
In the case $n=1$ it is convenient to define a smooth function $u:X\to \bbr$ to be   atomic if $\oper{log} |u|\in\loc (X)$, (c.f. [HS]).
\newline
\newline

If $u$ is atomic then the zero set $Z$ has measure zero in $X$ (see [HS]) 
so that the $\loc (X)$ extensions are unique. In particular,  the smooth
form $u^*(\th)$ on $X\sim Z$ has a unique $\loc (X)$ extension across $Z$,
and therefore defines a current on $X$.

\subheading{Definition 2.2} Let $u: X\to \bbr^n$ be an atomic function.
The {\bf divisor} of $u$ is the degree $n$ current $\Div (u)$ on $X$ defined by $$ 
\oper{Div} (u)\,\, := \,\, d( u^*\th ). 
$$
Atomicity is a weak condition which ensures the existence of a zero divisor.
Harvey and Semmes proved that a large class of smooth functions are atomic.
More specifically those functions which vanish algebraically and whose 
zero sets are not too big in the sense of Minkowski content are atomic. In 
particular real analytic functions whose zero sets have codimension $n$ are atomic.

\proclaim{Lemma 2.3}
Let $g$ be a smooth $GL (n,\bbr )$--valued function on an \od manifold $X$
and let $u: X\to \bbr^n$ be atomic. Then $v := ug$ is atomic and
$$
\oper{Div} (v)\,\, = \,\, \pm 1\, \oper{Div} (u),
$$
where $\pm 1 := \oper{sgn}\oper{det} (g)$ is constant on connected components of $X$.

\endproclaim

This result of [HS] allows one to extend the notion of divisor to sections of vector bundles. First, a  section $\nu$ of a smooth vector bundle $V\to X$ is called {\bf atomic} if
for each choice of local frame $e$ for $V$ the function $v$, defined by $\nu = ve$, is atomic.

\subheading{Definition 2.4}
Let $\nu$ be an atomic section of a rank $n$ bundle  $V\to X$.
 The {\bf  divisor,} $\Div (\nu )$, of $\nu$ is the
  $\OV$--twisted current on $X$ defined locally on an open subset $U$ of $X$ as follows. Choose  a local frame $e$ for $V$ over $U$ and let $v: U\to \bbr^n$ be the coordinate expression for $\nu$ determined by $e$. Then
$$ 
\Div (\nu ) \,\, := \,\, [e]\otimes \Div (v) \qquad\text{on } U.
$$
In particular, if $V$ is oriented, then $\Div (\nu )$ is a current on $X$.

As described in the Appendix the {\bf locally integrally flat} currents are those currents which can can expressed as  $R + dS$ where $R$ and $S$ are locally rectifiable. Furthermore the complex ${\wt{\SF}}^*_\LOC (X)$ of $\OV$--twisted currents on $X$ which are locally integrally flat may be used to compute the cohomology, $H^*(X, \wt\bbz_V )$, of $X$ with integer coefficients
twisted by $\SO_V$, i.e. $\wt\bbz_V := \SO_V\otimes_{\bbz_2}\bbz$.

\proclaim{Theorem 2.5}
 Let $\nu$ be an atomic section  of a real rank $n$ vector bundle $V\to X$. The  zero divisor, $\oper{Div} (\nu )\in \wt{\SF}^n_\LOC (X)$, of $\nu$  is an
 $\OV$--twisted  $d$--closed locally integrally flat current of degree $n$
 on $X$, whose support is contained in  the zero set of the section $\nu$. Furthermore, if $\mu$ is another atomic section of $V$, then there is an 
$\OV$--twisted locally rectifiable current $R$ so that
$$ 
\oper{Div} (\nu ) \,\, - \,\, \oper{Div} (\mu ) \,\, = \,\, dR.\tag 2.6
$$
That is, the cohomology class of $\oper{Div} (\nu )$ in $H^n (X , \wt\bbz_V )$ is well defined independent of the choice of section $\nu$. This class is the   twisted Euler class $\wt e$ of $V$. In particular, if $V$ is oriented, the 
Euler class $e\in H^n(X,\bbz )$ of $V$ is the cohomology
class of $\Div (\nu )$.
\endproclaim

\proclaim{Corollary 2.7}
Suppose that $\nu$ is an atomic section of an {\bf odd} rank bundle $V\to X$.
Then there is an $\OV$--twisted locally rectifiable current $R$ on $X$ so that 
$$ 
2 \, \Div (\nu ) \,\, = \,\, dR.
$$
Consequently the cohomology class $\wt e \in H^n (X , \wt\bbz_V )$ of $\Div (\nu )$ is a torsion class  
of order 2, (c.f. [MS]).
\endproclaim

\demo{Proof of Corollary}
Since the antipodal map on the even dimensional sphere $S^{n-1}$ is orientation
reversing, $\Div (-\nu ) = -\,\Div (\nu )$. The result now follows by applying equation (2.6).
\qed\enddemo

The proof of Theorem 2.5 relies on the following general remark.

\subheading{Remark 2.8}
 Let \pxtx denote the double cover $\SO_V\to X$. Note that  the pullback bundle $\wt V = p^* V$ on $\Xt$ is \od by choosing the orientation on $\Vt_{\wt x}$ to be the one determined by the point $\wt x \in\Xt$.
A current $T$ on $\wt X$ is {\bf odd} if $a_* T = \, -T$, where $a: \Xt\to\Xt$ is the natural involution. Then odd  currents on $\Xt$ are in 1--1 correspondence with $\OV$--twisted currents on $X$. 

\demo{Proof of Theorem 2.5}
The proof in the case that $V$ is oriented is given in [HS,5.1]. In the case that $V$ is nonorientable
let \pxtx be the double cover of Remark 2.8 and let  $\wt\nu$ denote the pullback of the section $\nu$ to $\wt V\to \Xt$. Since $\wt V$ is oriented the
theorem holds for the section $\wt\nu$. Furthermore, by Lemma 2.3, $\Div (\wt\nu )$ is an odd  current on $\Xt$ which corresponds to the $\OV$--twisted current $\Div (\nu )$ on $X$. Incorporating the double cover \pxtx into the proof of [HS, 5.1] we observe that the locally rectifiable current $\wt R$ on $\Xt$ can be chosen to be odd. 
Finally we let $R$ be the corresponding $\OV$--twisted current on $X$.
\qed\enddemo

The following structure theorem for  divisors, which is a corollary of
[HS,4.3], is proved in [Z1]. Set $Z := \oper{Zero} (\nu )$ and let
$$
\oper{Reg} Z \,\, = \,\ \{ x\in X \,\, : \,\, Z \,\,\text{is a codimension--} n\,\,  C^1 \text{ submanifold near } x\}
$$
denote the set of regular points of $Z$ and let $\oper{Sing} Z := Z \sim\oper{Reg} Z$ denote the set of singular points. Let $\{ Z_j \}$ denote the family of connected components of $\oper{Reg} Z$.

\proclaim{Theorem 2.9}
Let $\nu$ be an atomic section of $V\to X$. If $Z_j \subset \oper{spt} (\Div (\nu ))$, then 
$$ 
\SO_{TZ_j} \,\, \cong \,\, \OV\otimes{\SO_X} \bigm|_{Z_j},
$$
and, given such an isomorphism, the submanifold $Z_j$ defines an $\OV$--twisted current $[Z_j]$ by integration. Furthermore, there are integers $n_j\in\bbz$ such that
$$
\Div (\nu ) \,\, = \,\, \sum n_j\, [Z_j]\tag 2.10
$$
 as $\OV$--twisted currents on $X\sim \oper{Sing} Z$.
\endproclaim

\subheading{Remark 2.11}
By [HS, 4.3], the integers $n_j$ in equation (2.10) can be calculated as follows. Let $x\in Z_j$ and let $\SU$ be an open neighbourhood of $x$ in $X$. Choose orientations for $V$ and $TX$ over $\SU$. These orientations induce  orientations on $TZ_j$ and on the normal bundle $NZ_j$ over $\SU$. Let $\rho : S(V)\to S(\bbr^n )$ be an orientation preserving trivialization of the sphere bundle $S(V)$ over $\SU$. Then, for
almost all $x\in Z_j$, 
$$
n_j \,\, = \,\, \oper{Deg} \( \rho\circ\nu : S(N_xZ_j) \to S(\bbr^n )\)
$$
is the degree of the induced map between oriented $(n-1)$--spheres.

\subheading{B. Mod 2 divisors}
\newline
\newline

In this subsection the mod 2 divisor of an atomic section is defined to be the mod 2 reduction of the divisor of the section. The mod 2 reduction of an $\OV$--twisted locally integrally flat current is a mod 2 locally integrally flat current. At the cohomology level mod 2 reduction is simply the natural mapping $H^*(X, \wt\bbz_V)\to H^*(X,\bbz_2)$. The idea behind the definition
of a mod 2 current is to completely ignore orientation issues by declaring a current $T$ and its negative $-T$ to be the same. 
 Although they encode less information than  their twisted counterparts, mod 2 currents have the advantage that they pushforward under proper smooth maps. This fact will be particularly useful in Section 3.

\subheading{Definition 2.12}
Let $\SF^p_\LOC (X)$ denote the space of locally integrally flat currents of 
degree $p$ on $X$. 
Then the  space $\SF^{\oper{mod} 2}_\LOC (X)$ of {\bf mod 2 locally integrally flat currents} of degree $p$ on $X$ is defined to be the quotient $\SF^p_\LOC (X) / 2\SF^p_\LOC (X)$. The natural mapping $\SF_\LOC (X) \to \SF^{\oper{mod} 2}_\LOC (X)$ is called {\bf mod 2 reduction}.
\newline
\newline

The spaces $\SF_\LOC (X)$ and $\SF^{\oper{mod} 2}_\LOC (X)$ are studied in the Appendix. In particular we show there that the complex $\SF^{\oper{mod} 2}_\LOC (X)$ may be used to compute the cohomology, $H^*(X,\bbz_2)$, of $X$ with $\bbz_2$ coefficients. In Subsection A we saw that $\Div (\nu )\in \wt{\SF}^n_\LOC (X)$ is an $\SO_V$--twisted locally integrally flat current. Now, by Lemma A.13, there is a canonical isomorphism
$$
\wt \SF_\LOC (X) \bigm/ 2 \wt\SF_\LOC (X) \,\, \cong \,\, \SF^{\oper{mod} 2}_\LOC (X). \tag 2.13
$$
The induced mapping $\wt\SF_\LOC (X)\to\SF^{\oper{mod} 2}_\LOC (X)$ is also called {\bf mod 2 reduction}.

\subheading{Definition 2.14}
The {\bf mod 2 divisor}, $\Div^{\oper{mod} 2} (\nu )\in \SF^{\oper{mod} 2}_\LOC (X)$, of an atomic section $\nu$ 
of $V\to X$ is defined to be the mod 2 reduction of the $\OV$--twisted
current $\Div (\nu )\in \wt\SF_\LOC (X)$. 
\newline
\newline

Now the mod 2 version of Theorem 2.5 is immediate.

\proclaim{Theorem 2.15}
 Let $\nu$ be an atomic section  of a real rank $n$ vector bundle $V\to X$. 
The   mod 2  divisor, $\oper{Div}^{\oper{mod 2}} (\nu) $,  is a
   $d$--closed  mod 2 locally integrally flat current of degree $n$
 on $X$, whose support is contained in  the zero set of the section $\nu$.
 Furthermore, if $\mu$ is another atomic section of $V$, then there is a 
mod 2 locally rectifiable current $R$ so that
$$ 
\oper{Div}^{\oper{mod 2}} (\nu ) \,\, - \,\, \oper{Div}^{\oper{mod 2}}
 (\mu ) \,\, = \,\, dR.
$$
That is, the cohomology class of $\oper{Div}^{\oper{mod 2}} (\nu )$ in $H^n (X , \bbz_2 )$ is well defined independent of the choice of section $\nu$.
 This class is the   mod 2  Euler (or top Stiefel-Whitney) 
 class, $w_n$, of $V$. 
\endproclaim

\subheading{Remark 2.16}
By definition the divisor $\Div (\nu)\in \wt\SF_\LOC (X)$ is determined by a 
collection of local divisors $\Div (v_\a )$ defined on open subsets $U_\a$ of $X$. These local divisors satisfy $\Div (v_\a ) = \pm\Div (v_\b )$ on $U_\a\cap U_\b$. The mod 2 divisor $\Div^{\oper{mod} 2} (\nu)$ is the mod 2 current which is naturally associated to this collection of local divisors.
\newline
\newline

\subheading{Example 2.17}
In general a current representative  for a mod 2 divisor  is not $d$--closed.
 In fact it is easy to construct sections $\nu$ for which there are {\it no} $d$--closed {\it current} representatives
 $T\in \SF_\LOC (X)$ of 
$\Div^{\oper{mod} 2} (\nu)$, (though of course $\Div (\nu)$ is a $d$-closed $\OV$-{\it twisted} current representative of $\Div^{\oper{mod} 2} (\nu)$). This can be done as follows. 

Let $\nu$ be an atomic section of a nontrivial real line bundle $L\to X$. 
Choose a metric on $L$. Then the divisor of $\nu$ is the $d$--closed $\SO_L$--twisted current $\Div (\nu) = d (\tfrac 12 \f{\nu}{|\nu |})
\in \wt{\SF}^1_\LOC (X)$. The mod 2 reduction of $\Div (\nu)$ is the mod 2 divisor $\Div^{\oper{mod} 2} (\nu)$ which represents the first Stiefel--Whitney
class $w_1 (L)\in H^1(X, \bbz_2)$. A current representative $T\in \SF^1_\LOC (X)$ of $\Div^{\oper{mod} 2} (\nu)$ can be constructed as follows. For simplicity
we assume that $0$ is a regular value of $\nu$. Let $Z = \oper{Zero} (\nu)$. 
Choose an auxillary   section $\mu$ of $L\to X$ so that $W = \oper{Zero} (\mu )$ is a submanifold which is transverse to $Z$. Let $p : \wt X \to X$ be the double cover
orienting $L$ and let $\wt\nu$, $\wt\mu : \Xt \to\bbr$ be the pullbacks of the sections $\nu$, $\mu$. Set $\chi := \f { \wt\mu}{|\wt\mu |}$. Now $\chi \Div (\wt\nu )\in \SF^1_\LOC (\Xt )$ is well defined and $T := \tfrac 12 p_* (\chi \Div (\wt\nu ))\in \SF^1_\LOC (X)$ is a current on $X$ whose mod 2 reduction is 
$\Div^{\oper{mod} 2} (\nu)$. Finally $dT = 2 [Z\cap W] \neq 0$ on $X$. 

Next suppose that $S = T - 2R$ is a current representative of $\Div^{\oper{mod} 2} (\nu)$ for which $dS = 0$. This forces $[Z\cap W] = dR$ to be zero in $H^2 (X,\bbz)$. But this is not always the possible,
since if $L\to \bbr\bbp^3$
 is the tautological line bundle then  $[Z\cap W]$ can be chosen to be the generator  $\bbr\bbp^1$ of $H^2(X,\bbz)\cong H_1(X,\bbz) = \bbz_2$.

\subheading{3. Linear dependency currents}
\newline
\newline

In this section we associate to each  atomic collection $\nu$ of $n-q+1$ sections of a real rank $n$ bundle $F\to X$ a degree $q$ current on $X$ which is supported on the linear dependency set 
of the collection of sections. This current will be called the {\bf linear dependency current} of the collection $\nu$. In all cases the linear dependency
current exists as  a mod 2 current. However, in the case where 
 $q$ is odd or where $q=n$, it can also be defined to be  an $\SO_F$--twisted
current. Note that if $q=n$ the linear dependency current is simply the divisor of the section $\nu$, (see Section 2).
 Henceforth we assume that $q<n$.

The linear dependency currents are defined using the  construction of such  currents in [HL3] which we now briefly recall.  Let $F\to X$ be a
real rank $n$ vector bundle. Fix $q\in \{ 1,2,...,n-1\} $ and let $ \nu = (\nu_1,...,\nu_m)$ be a collection of $m = n-q+1 > 1$
sections of $F\to X$. (Such collections will always be ordered.) These sections define a bundle map 
$\nu: \ur^m \to F$ by $\nu (t_1,...,t_m) := \sum\limits_{i=1}^m t_i\nu_i$
which drops rank on the set where $\nu_1,...,\nu_m$ are linearly dependent.
Let $\pi:\bbp (\ur^m)\to X$ denote the trivial bundle of $(m-1)$--dimensional
real projective spaces and let $U\subset \ur^m$ be the tautological line bundle
over $\bbp(\ur^m)$. Using $\pi$ to pull back the bundle map 
 $\nu : \ur^m\to F$ to $\bbp (\ur^m )$ and then restricting to the subbundle $U\subset \ur^m$ we obtain an {\bf induced section}  $\wt\nu$ of the bundle $H := \Hom (U,\pi^*F)$ over $\bbp(\ur^m) \equiv \bbp(\bbr^m)\times X$. By construction the projection by $\pi$ to $X$ of the zero set of $\wt\nu$ is the linear dependency set of $\nu_1,...,\nu_m$. 

\subheading{Definition 3.1} 
The collection $\nu$ of sections $\nu_1,...,\nu_m$ is called {\bf atomic} if the induced section $\wt\nu$ of $H\to\bbp(\ur^m)$ is atomic. 
\newline
\newline

If the collection $\nu$ is atomic, then the $\SO_H$--twisted divisor current,
$\Div (\wt\nu )$, and its mod 2 reduction, $\Div^{\oper{mod} 2} (\wt\nu )$,
 are well  defined.
The  linear dependency current is defined to be the current pushforward of 
 $\Div^{\oper{mod} 2} (\wt\nu )$, or whenever possible the pushforward  of
  $\Div (\wt\nu )$.

Generally speaking it is not possible to pushforward twisted currents on 
$\bbp(\bbr^m)\times X$ to $X$ by the projection $\pi$. However the pushforward
by $\pi$ of an ${\Cal O}_{\bbp(\bbr^m)}\otimes\pi^*\SO_F$--twisted current on
$\bbp(\bbr^m)\times X$ is well defined and is an $\SO_F$--twisted current on
 $X$. This observation together with the following elementary lemma will be used to determine
when it is possible to pushforward the ${\Cal O}_H$--twisted current $\Div(\wt\nu)$.

\proclaim{Lemma 3.2} If $n\equiv m$ mod 2 then there is a canonical
 isomorphism
$$
{\Cal O}_H  \,\, \cong \,\, \SO_{\bbp(\bbr^m)}\otimes\pi^*\SO_F.
$$
\endproclaim

\demo{Proof} 
 First recall that if $V$ and $W$ are oriented finite dimensional
 vector spaces then there is a canonical choice of orientation on $V\otimes W$.
Furthermore if the dimension of $W$ is even this choice is independent of the orientation on $V$ and so  there is a canonical isomorphism
 $\SO_W\cong\SO_{V\otimes W}$. 
The canonical choice of orientation on $V\otimes W$ is defined as follows. Choose ordered bases  $v=(v_1,...v_p)$ and $w=(w_1,...,w_q)$ for $V$ and $W$. 
Then the canonical orientation on $V\otimes W$ is given by the ordered
basis
$$
v\otimes w \,\,=\,\,
(v_1\otimes w_1,v_2\otimes w_1,...,v_p\otimes w_1,v_1\otimes w_2,...v_p\otimes w_2,...,v_p\otimes w_q).\tag 3.3
$$
Secondly recall [MS] that there is a canonical vector bundle isomorphism $U^*\otimes \UP \cong T\bbp(\bbr^m)$, where  $U^\bot$ denotes the orthogonal complement of $U$ in $\ur^m$.

So if $n$ and $m$ are both even, then the result is true since there are
canonical isomorphisms $\SO_H \cong \SO_{\pi^*F}$  (as $n$ is even) and $\SO_{\bbp(\bbr^m)}= \SO_{U^*\otimes\UP}
\cong \ur$. Here the last isomorphism is well defined by sending
 $[u^*\otimes u^\bot]$ to 1, where $u^\bot$ is chosen so that $(u,u^\bot)$ 
is a positively oriented frame for $U\oplus \UP =\ur^m$. Similarly if $n$
 and $m$ are both odd then $\SO_H \cong \SO_{U^*}\otimes\SO_{\pi^*F}$ and 
since $m-1$ is even 
$\SO_{U^*}\cong \SO_{U^*\otimes\UP} = \SO_{\bbp(\bbr^m)}$, as required.
\qed\enddemo

\subheading{A. Linear dependency currents (q odd)}
\newline
\newline

 Throughout this subsection we assume that $q<n$ is odd  and hence $m\equiv n$ (mod 2). 
Then, by   Lemma 3.2, the pushforward by $\pi$ of the $\SO_{\bbp(\bbr^m)}
\otimes\pi^*\SO_F$-twisted current $\Div (\wt\nu )$ on $\bbp(\bbr^m)\times X$
 exists and is an $\SO_F$-twisted current on $X$.

\subheading{Definition 3.4}
Let $q$ be odd. The 
 {\bf  linear dependency current}, $\BL\BD (\nu)$, of  an atomic collection $\nu$ of $n-q+1$ sections of $F\to X$ is the $\SO_F$-twisted current on $X$ defined by
$$
\BL\BD (\nu ) := \pi_* (\oper{Div} (\widetilde\nu )).
$$

\subheading{Remark 3.5}
The following equivalent definition of the  linear dependency
current is often useful, especially when $n$ and $m$ are both odd as in this case the fibres  $\Hom (U, F_x )$ and $\bbp (\ur^m )_x$ are nonorientable.
 Let $p: S(\ur^m )\to\bbp (\ur^m)$ be the double cover 
by the unit sphere and let $\rho := \pi\circ p : S(\ur^m)\to X$. Let $NS(\ur^m)$ denote the normal bundle to $S(\ur^m)$ in $\ur^m$, with its canonical orientation. Then, as above, there is associated to the collection $\nu$ a section $\widehat\nu$ of the bundle $\widehat H := \Hom (NS(\ur^m),\rho^*F)$ over $S(\ur^m)$. Note that the sections $\widehat\nu$ and $\nut$ are simultaneously atomic. As above there is a canonical isomorphism $\SO_{\widehat H}\cong\rho^*\SO_F$. Consequently $\Div (\widehat\nu )$ is a well defined $\rho^*\SO_F$--twisted current on $S(\ur^m )$. Then, if $q$ is odd, 
we have that
$$
\BL\BD (\nu ) = \tfrac 12\, \rho_* (\oper{Div} (\widehat\nu )).\tag 3.6
$$

We verify (3.6) as follows. First note that $\widehat H = p^*H$ and that $\widehat\nu = p^*\nut$. Since $n\equiv m$ (mod 2),  the $\rho^*\SO_F$--twisted current
$\Div (p^* \nut )$ is even in that
$$
a_* \Div (p^* \nut )\,\, = \,\, \Div (p^* \nut )\qquad\text{on } S(\ur^m),
$$
where $a : S(\ur^m)\to S(\ur^m)$ is the antipodal map. Now even currents on $S(\ur^m)$ are in 1--1 correspondence with $\SO_{\bbp(\bbr^m)}$--twisted
 currents on $\bbp (\ur^m )$. In particular, if $q$ is odd, 
$$
\Div (\nut )\,\, = \,\, \tfrac 12 \, p_* (\Div (p^* \nut ))\qquad\text{on } \bbp (\ur^m).
$$
This fact immediately implies (3.6). 
Note that, if $q$ is even, then $\Div (p^* \nut )$ is an odd twisted current
on $S(\ur^m)$ and so its current pushforward is zero. In general odd currents on $S(\bbr^m)$ are in 1--1 correspondence with $\SO_U\otimes \SO_{\bbp(\bbr^m)}
$--twisted currents on $\bbp (\bbr^m )$, see [Z2]. 
\newline
\newline

The following result generalizes Theorem 2.5.

\proclaim{Theorem 3.7} 
Let $F\to X$ be a real rank $n$ bundle  and let $q$ be odd. For each atomic collection $\nu$ of $n-q+1$ sections  of $F\to X$  the linear dependency current $\BL\BD (\nu )$ is an  $\orF$--twisted $d$--closed
locally integrally flat current of degree $q$ on $X$ whose support is contained in the linear dependency set of the collection of sections. Furthermore, if
 $\mu$ is another atomic collection of sections of $F\to X$, then there
is an $\orF$--twisted  locally rectifiable current  $R$ so that
$$ \BL\BD (\nu ) \,\, - \,\, \BL\BD (\mu ) \,\, = \,\, dR.  $$
That is, the cohomology class of $\BL\BD (\nu )$ in $H^q(X,\wt\bbz_F )$ is well
defined independent of the choice of sections.
\endproclaim

\demo{Proof}
Since the pushforward of a locally rectifiable current is locally rectifiable
the current $\BL\BD (\nu )$ inherits its properties from those of the
divisor of the  induced section $\wt\nu$ (see Theorem 2.5). 
\qed\enddemo

\subheading{Note}
In the next section the cohomology class of $\BL\BD (\nu )$ is shown to be 
$\wt W_q\in H^q(X,\wt\bbz_F)$, the (twisted) integer Stiefel--Whitney class of $F$, whose mod 2 reduction is the standard Stiefel--Whitney class $w_q\in H^q(X,\bbz_2)$ of $F$.

\proclaim{Proposition 3.8}
Let $\nu : \ur^m\to F$ be as above and let $\psi: \ur^m\to\ur^m$ and $\varphi : F\to F$ be bundle isomorphisms. Then the collection of sections corresponding to the bundle map $\varphi\circ\nu\circ\psi : \ur^m\to F$ is also atomic. Furthermore,
if $q$ is odd, 
$$
\BL\BD (\varphi\circ\nu\circ\psi ) \,\, = \,\, \oper{sgn}\oper{det} (\psi ) \, 
\oper{sgn}\oper{det} (\varphi ) \, \BL\BD (\nu ),
$$
as $\SO_F$--twisted currents on $X$. 
\endproclaim

The proof of the proposition will be given at the end of this subsection.
The following result generalizes Corollary 2.7.

\proclaim{Corollary 3.9}
Under the same hypothesis as in Theorem 3.7, with $q < n$,   there is an $\SO_F$--twisted locally rectifiable
current $R$ on $X$ so that 
$$
2\, \BL\BD (\nu ) \,\, = \,\, dR. 
$$
Consequently the cohomology class $\wt W_q$ of  $\BL\BD (\nu )$ in $H^q(X,  \wt\bbz_F )$ is  a torsion class of order 2.
\endproclaim

\demo{Proof of Corollary}
Define $\psi : \ur^m\to\ur^m$ by $\psi (t_1,...,t_m) \, = \, (-t_1, t_2,...,t_m)
$ and let $\mu := \nu \circ \psi$. Then, by Proposition 3.8, $\BL\BD (\mu ) = \, - \,\BL\BD (\nu )$. The result now follows from Theorem 3.7. For an alternate proof see Theorem 4.10.
\qed\enddemo

Next we study the case $q=1$ in more detail.

\subheading{Remark 3.10}
Let $\SO_F^\bbr := \SO_F\otimes_{\bbz_2} \bbr$ denote the orientation line
bundle of $F$. The divisor $\Div (\s) = d(\tfrac 12 \frac {\s}{|\s |})$ of an atomic section $\s$ of $\SO_F^\bbr$ is called the {\bf orientation current} of $F$ associated with $\s$. Note that $\Div (\s )$ is an $\SO_F$--twisted current  whose  cohomology class in $H^1 (X,\wt \bbz_F)$ is $\wt W_1$.  
\newline
\newline

To each collection $\nu$ of $n$ sections of a rank $n$ bundle $F$ there is an associated section $\s$ of $\SO_F^\bbr$ well defined as follows. Choose a local frame $f$ for $F$ and let $A$ be the matrix defined by $\nu = A f$. Then
$\s := [f]\otimes \oper{det} A$.

\proclaim{Proposition 3.11. The case q = 1}
Let $\nu$ be an atomic collection of $n$ sections of a rank $n$ bundle $F$.
Suppose that the associated section $\s$ of $\SO_F^\bbr$ is also atomic.
Then
$$
\BL\BD (\nu ) \,\, = \,\, \Div (\s ) $$
as $\SO_F$--twisted currents on $X$. 
\endproclaim

\demo{Proof}
Choose a local frame $f$ for $F$ and define
$A$ by $\nu = A f$.  Then the local
expression for the induced section $\widehat\nu$ of $\rho : \oper{Hom} (NS(\ur^n), \rho^*F)\to S(\ur^n)$ is the mapping $\psi : X\times S(\bbr^n) \to \bbr^n$ defined by $\psi (x,y) = y A(x)$. Let $\theta$ denote the normalized solid angle kernel on $\bbr^n$. 
Now, by the Change of Variables and Stokes' Theorems, 
$$
\int\limits_{\rho^{-1} (x)} \psi^*\theta \,\, = \,\, \f{\oper{det} A(x) }{| \oper{det} A(x) |}
\qquad\text{for each } x\notin \oper{Zero} (\oper{det} A).
$$    
This implies the  Proposition  since  $\Div(\s) = [f]\otimes d \(\tfrac 12 \f{\oper{det} A}{| \oper{det} A|}\)$ and $\BL\BD (\nu) = [f]\otimes  \tfrac 12\rho_* (d\psi^*\theta ))$. 
\newline
\qed\enddemo

\demo{Proof of Proposition 3.8} By (3.6) it suffices to show that
$$
\Div (\widehat{\varphi\circ\nu}) \,\, = \oper{sgn}\oper{det} (\varphi ) \, \Div (\widehat\nu )\qquad\text{on } S(\ur^m)\tag 3.12
$$
 and, if $q$ is odd, that
$$
\Psi_*\Div(\widehat{\nu\circ\psi})\,\,=\,\, \oper{sgn}\oper{det} (\psi) \, \Div (\widehat\nu ),\qquad\text{on } S(\ur^m), \tag 3.13 
$$
where $\Psi : S(\ur^m) \to S(\ur^m)$ is the diffeomorphism induced by $\psi$. 
Now since $\widehat{\varphi\circ\nu} = \varphi\circ\widehat\nu$, (3.12) follows from Lemma 2.3. To prove (3.13) let $\mu=\nu\circ\psi$ and note that
the pullback of the section $\widehat\nu$ of $\widehat H$ by $\Psi$ is a section $\Psi^*\widehat\nu$ of $\Psi^* \widehat H = \oper{Hom} (\Psi^*NS(\ur^m),\rho^*F)$. Let $\psi^*: \oper{Hom}  (\Psi^*NS(\ur^m),\rho^*F) \to \oper{Hom} (NS(\ur^m),\rho^*F)$ be the bundle isomorphism defined by $\psi^* (\a ) := \a\circ\psi$. Then
$$
\widehat\mu \,\, = \,\, \psi^*\, ( \Psi^*\widehat\nu ).
$$
Clearly $\widehat\nu$, $\Psi^*\widehat\nu$ and $\widehat\mu$ are simultaneously atomic. 
 Since   equation (3.13)   is local on $X$ we can assume that $F$ and $X$ are oriented. Then $\widehat\mu$ and $\widehat\nu$ are sections of the oriented bundle $\widehat H$ over the oriented manifold $S(\ur^m)$. Furthermore
since $q$ is odd  the orientation induced on $\Psi^*\widehat H$ by the diffeomorphism $\Psi$ is the same as that induced by the bundle isomorphism $\psi^*$.
Equation (3.13) now follows immediately by applying Lemma 2.3 to $\psi$, the
Change of Variables Theorem to $\Psi$, and noting that 
 $\oper{sgn}\oper{det} (D\Psi) = \oper{sgn}\oper{det} (\psi)$.
\qed\enddemo

\subheading{B. Mod 2 linear dependency currents}
\newline
\newline

Since  mod 2 currents can always be pushed forward by proper maps we can
define the mod 2 linear dependency current for an even as well as an odd number of sections.

\subheading{Definition 3.14}
Suppose $1\leq q\leq n$.  The {\bf mod 2 linear dependency current},
$\BL\BD^{\oper{mod} 2} (\nu )$, of an atomic collection $\nu$ of $m=n-q+1$ sections  of $F\to X$ 
is defined to be the current pushforward of the mod 2 divisor of the induced section $\nut$ of $H\to \bbp(\ur^m )$,
$$
\BL\BD^{\oper{mod} 2} (\nu ) \,\, := \,\, \pi_* (\Div^{\oper{mod} 2} (\nut )).
$$

\subheading{Note}
If $q$ is odd or $q=n$ the mod 2 linear dependency current, $\BL\BD^{\oper{mod} 2} (\nu )$, is the mod 2 reduction of the $\orF$--twisted linear dependency
current, $\BL\BD (\nu )$.
\newline
\newline

The mod 2 analogues of Theorem 3.7 and Proposition 3.8 hold. In particular

\proclaim{Theorem 3.15} 
 For each atomic collection $\nu$ of $n-q+1$ sections  of a real
rank $n$ bundle  $F\to X$  the linear dependency current
 $\BL\BD^{\oper{ mod 2}} (\nu )$
 is a $d$--closed mod 2
locally integrally flat current of degree $q$ on $X$ whose support is contained in the linear dependency set of the collection of sections. Furthermore, if
 $\mu$ is another atomic collection of sections of $F\to X$, then there
is a mod 2  locally rectifiable current  $R$ so that
$$ \BL\BD (\nu )^{\oper{ mod 2}} \,\, - \,\, \BL\BD^{\oper{ mod 2}} (\mu ) \,\, = \,\, dR.  $$
That is, the cohomology class of $\BL\BD^{\oper{ mod 2}} (\nu )$
 in $H^q(X,\bbz_2 )$ is well
defined independent of the choice of sections.
\endproclaim

\subheading{Note}
 In the next section the cohomology class of  $\BL\BD^{\oper{mod} 2} (\nu )$ is shown to be $w_q(F)\in H^q (X,\bbz_2)$, the $q$th Steifel--Whitney class of $F$.

\subheading{C. The structure of linear dependency currents}
\newline
\newline

The following result concerning the structure of the twisted linear dependency
current builds on Proposition 2.8 of [HL] and Theorem 2.9 above. Let $q$ be odd and let $\nu$ be an atomic collection of $m = n-q+1$ sections of $F\to X$. Suppose that the zero set,
$Z(\nut )$, of  $\nut$ is a smooth submanifold of $\bbp (\bbr^m) \times X$. Let $\{ Z_j \} $ denote the connected components
of $Z(\nut )$. Then, by Theorem 2.9, there are integers $n_j\in\bbz$ so that
$$
\Div (\nut ) \,\, = \,\, \sum n_j [Z_j]\tag 3.16
$$
as $\SO_{\bbp(\bbr^m)}\otimes\pi^*\SO_F$--twisted currents on $\bbp(\ur^m)$. By [HL, 2.8],
the subset
$$
RK_{m-1} (\nu ) \,\, := \,\, \{ x\in X \,\, : \,\, \oper{rank} \nu_x \,\, = \,\, m-1 \}
$$
of the linear dependency set of $\nu$ is a locally rectifiable set. Let 
$RK_j := RK_{m-1} (\nu )\cap \pi (Z_j)$. If $n_j\neq 0$ then
$$
\SO_{TRK_j}\,\, \cong\,\, \SO_F\otimes{\SO_X}\bigm|_{RK_j}
$$
(wherever it makes sense) and, given an isomorphism of these two bundles,
$RK_j$ defines an $\SO_F$--twisted current $[RK_j]$ by integration.
Arguing as in the proof of [HL, 2.8] it follows that $\pi_* [Z_j] = [RK_j]$. Consequently we have the following

\proclaim{Proposition 3.17}
Let $\nu$ be as above. Then 
$$
\BL\BD (\nu ) \,\, = \,\, \sum n_j [RK_j]
$$
as $\SO_F$--twisted currents on $X$, where the integers $n_j$ are
given by (3.16).
\endproclaim

Next we examine the structure of the twisted and mod 2 linear dependency
currents in the case that $m-1$ of the $m$ sections are everywhere linearly independent. 

\proclaim{Theorem 3.18}
Let $\mu_1,...,\mu_m$ be a collection of $m=n-q+1$ sections of $F\to X$. Suppose that
$\mu_1,...,\mu_{m-1}$ are everywhere linearly independent. Choose a metric on $F$ and let $\mu_m^\bot$ denote the  projection of $\mu_m$ onto the orthogonal 
complement of  $\mu_1,...,\mu_{m-1}$ in $F$. Then $\mu_m^\bot$ is atomic if and only if the induced section $\wt\mu$ of $H\to\bbp (\ur^m)$ is atomic. Furthermore, 
if $ q$  is odd, then
$$ \BL\BD (\mu )\,\, = \,\,\Div (\mu_m^\bot ) \qquad\text{as } \SO_F \text{--twisted currents on } X, \tag 3.19
$$
and, 
for any $q$,
$$ \BL\BD^{\oper{mod} 2} (\mu )\,\, = \,\,\Div^{\oper{mod} 2} (\mu_m^\bot ) \qquad\text{as mod 2 currents on } X. \tag 3.20
$$

\endproclaim

\demo{Proof}
We present the proof in the twisted case. The mod 2 case follows similarly.
 The first step is to choose local coordinates and frames
and to relate the local coordinate expression for the induced section $\wt\mu$
of $H\to X\times\bbp (\bbr^m)$ to that of the section $\mu_m^\bot$.

First note that the linear dependency set of $\mu_1,...\mu_m$ is equal to the zero set $Z(\mu_m^\bot )$ of $\mu_m^\bot$. Fix a point $x_0\in Z(\mu_m^\bot )$
and let $\SU$ be a sufficiently small open neighbourhood of $x_0$ in $X$. Choosing orientations for $TX$ and $F$ over $\SU$ we can regard $\BL\BD (\mu )$ and $\Div (\mu_m^\bot )$ as currents on $\SU$. Now, since $\mu_1,...,\mu_{m-1}$
are linearly independent, there is precisely one point $\wt x$ of the zero set 
$Z(\wt\mu )$ of $\wt\mu$ in $X\times\bbp (\bbr^m )$ lying over each point $x$ of $Z(\mu_m^\bot )$ in $X$. Let $\SW \subset \pi^{-1} (\SU )$ be a sufficiently small open neighbourhood of $\wt {x_0}$ in $X\times\bbp (\bbr^m )$ which 
contains $Z(\wt\mu ) \cap \pi^{-1} (\SU )$. 

Choose the coordinate chart $\psi : \bbr^{m-1} \to \bbp (\bbr^m )$ defined by $\psi (s) = [s,1]$ and the local frame $u$ for $U\to\bbp (\bbr^m )$ over $\psi (\bbr^{m-1})$
defined by $u(s) = (s,1) \in U_{[s,1]}\subset \ur^m$. Note that the orientations induced on $T\bbp(\bbr^m)$ and $U$ by $\psi$ and $u$ are compatible and that $\Cal W \subset \psi(\bbr^{m-1})\times \Cal U$ since $\mu_1,...,\mu_{m-1}$ are linearly independent on $\Cal U$.  Also note that the orientations on $U$ and $F$ induce a natural orientation on $H$ over $\SW$. 

Choose a positively oriented local frame $f_1,...,f_n$ for $F$ over $\Cal U$  so that $f_i = \mu_i$ for $1\leq i\leq  m-1$ and $f_i\bot \oper{Span} \{ f_1,...,f_{m-1} \}$ for $i > m-1$. Define $a_i :\Cal U \to \bbr$ by 
$$
\mu_m = \sum\limits_{i=1}^n a_i f_i.
$$
Let $a' = (a_1,...,a_{m-1} )$ and  $a'' = (a_m,...,a_n)$.

Then the local coordinate expression for $\mu_m^\bot$ in terms of the local frame $f_m,...,f_n$ is $a'' : \SU\to \bbr^q$ and the coordinate expression for
$\wt\mu$ in terms of the frames $u,f$ is 
  $(a'',s + a') : \SU\times\bbr^{m-1}\to \bbr^n$. These two coordinate expressions can be related as follows. 
 Let $\Psi: \SU\times\bbr^{m-1} \to \SU\times\bbr^{m-1}$ be the orientation preserving change of variables
$\Psi (x,s) = (x,s+a'(x))$ and let $\oper{Id} : \bbr^{m-1}\to\bbr^{m-1}$ denote the identity map.  Then 
$$(a'',s + a') \,\, = \,\, (a''\times\oper{Id})\circ \Psi.$$ 
By the Change of Variables Theorem and Lemma 3.21 below it follows that $\pi_* (\Div (\wt\mu )) = \Div (\mu_m^{\bot} )$ as required.
\qed\enddemo

The following elementary fact about divisors is
 included for the sake of completeness.

\proclaim{Lemma 3.21}
Let $X$ be an oriented manifold and $f: X\to \bbr^n$ a smooth map.
 Let $\oper{Id} : \bbr^m\to\bbr^m$ be the identity map and let $\pi : X\times \bbr^m \to X$ denote projection onto $X$. Then $f$ is atomic if and only
if $f\times \oper{Id}$ is atomic. Furthermore, if $X\times \bbr^m$ is given 
the induced orientation, then
$$
\pi_* (\Div (f\times \oper{Id})) \,\, = \,\, \Div (f).
$$
\endproclaim

\demo{Proof}
 By induction we may assume that $m=1$. Let $t$ denote the coordinate
 on $\bbr$. 
First, $\f{dtdf^I}{|(t,f)|^{p+1}}\in \loc(X\times\bbr )$ for $p=|I|\leq n-1$ 
iff $\log |v|\in\loc (X)$ in case $p=0$, and iff $\f{df^I}{|f|^p}\in\loc (X)$ in case $p>0$. This is because
$$\int\limits_{|t| < R} \f{dt}{(|t| + |f|)^{p+1}} \,\, = \,\, \f{2}{|f|^p}\,\int\limits_0^{R/|f|} \f{ds}{(s+1)^{p+1}}.
$$
Of course $\f{df^I}{|(t,f)|^p}$ is dominated by $\f{df^I}{|f|^p}$.
Consequently $f\times\oper{Id}$ is atomic if and only if $f$ is atomic.

Let $\th_n$ denote the normalized solid angle kernel on $\bbr^n$, and recall that $\Div (f) = d(f^* \th_n )$ where $d$ is  exterior differentiation of generalized forms. Let $\p S$ denote the boundary of a current $S$. Then $dS = (-1)^{k+1} \p S$, where $k = \oper{deg} S$. Consequently  $\pi_*d \, = \, - d\pi_*$.
Therefore it suffices to show that
$$
\pi_* \( (f\times\oper{Id})^*\th_{n+1} \) \,\, = \,\, -f^*\th_n.
$$
Let $\om_n = \oper{Vol} (S^{n-1})$ and  $\lambda (f) = df_1\wedge...\wedge df_n$.  Then
$$
(f\times\oper{Id})^*\th_{n+1} \,\, = \,\, \f{1}{\om_n} \, \f{t\, \lambda (f)}{(t^2 + |f|^2)^{n/2}} \,\, - \,\, \f{\om_{n-1}}{\om_n}\, \f{|f|^{n-1} dt}
{(t^2 + |f|^2)^{n/2}} \, f^*\th_n.
$$
Now, since the pushforward of $(f\times\oper{Id})^*\th_{n+1}$ by $\pi$ is equal
to the integral of $(f\times\oper{Id})^*\th_{n+1}$ over the fibres of $\pi$,
$$
\align
\pi_*\((f\times\oper{Id})^*\th_{n+1}\) 
\,\,&=\,\, 
-\, \f{\om_{n-1}}{\om_n}\, \int\limits_{-\infty}^{\infty}\, \f{|f|^{n-1} dt}{(t^2 + |f|^2)^{n/2}}\, \, f^*\th_n \\
&= \,\,
 -\,\f{2\om_{n-1}}{\om_n}\, \int\limits_{0}^{{\pi}/2}\, \oper{cos}^{n-2} t\, dt
\, \, f^*\th_n \,\, = \,\,-\, f^*\th_n,\endalign
$$
as required.
\qed\enddemo

\subheading{4. Stiefel--Whitney currents}
\newline
\newline

The purpose of this section is to identify the cohomology class of a linear dependency current. First we consider the mod 2 case, and recall from Example 10
of the Appendix that $\bbz_2$ cohomology can be computed using mod 2 integrally flat currents,
$\SF_\LOC^{\oper{mod} 2} (X)$.

\proclaim{Theorem 4.1}
Given an atomic collection $\nu_1,...,\nu_m$ of $m = n-q+1$ sections of a real rank $n$ vector bundle $F$ over $X$ the mod 2 linear dependency current
$\BL\BD^{\oper{mod} 2} (\nu)\in \SF^{\oper{mod} 2}_\LOC (X)$ represents the $q$th Stiefel--Whitney class $w_q (F)\in H^q (X,\bbz_2 )$. 
\endproclaim

Before proving this result we note that the analogue of a theorem of Bott for complex vector bundles and Chern classes is valid for real vector bundles and 
Stiefel--Whitney classes. Let $w(E) := 1+ w_1(E) + ... + w_m (E)$ denote the total Stiefel--Whitney class of a real rank $m$ bundle $E$ over $X$. Let $U$
denote the universal line bundle on the projectivization $\bbp (E)$, and let $a := w_1 (U)\in H^1 (\bbp (E),\bbz_2)$ denote the first Stiefel--Whitney class of $U$ on $\bbp (E)$. Let $\pi : \bbp(E)\to X$ denote the natural projection, and let $\pi_* : H^{q+m-1} (\bbp (E),\bbz_2)\to H^q(X,\bbz_2)$ denote the induced map.

\proclaim{Lemma 4.2}
$$
\pi_* \( (1+a)^{-1}\) \,\, = \,\, w(E)^{-1}.
$$
\endproclaim

\demo{Proof}
Choose an inner product for $E$. Let $\BE$ denote the pullback of the bundle 
$E$ to $\bbp (E)$, and let $U^\bot$ denote the orthogonal bundle to $U\subset \BE$. The  product formula for  Stiefel--Whitney classes implies that
$$
w(\BE ) \,\, = \,\, w(U)\, w(U^\bot )\tag 4.3
$$
so that 
$$w(U)^{-1} \,\, = \,\, w(\BE )^{-1} \, w(U^\bot ).\tag 4.4$$
Since the fibre dimension of $\bbp (E)$ is $m-1$,
$$\pi_* (w_j (U^\bot )) \,\, = \,\, 0 \qquad\text{if } j < m-1.\tag 4.5
$$
Therefore, $\pi_*(w(U^\bot)) = \pi_* (w_{m-1} (U^\bot))$.
It remains to show that 
$$
\pi_* (w_{m-1} (U^\bot)) \,\, = \,\, 1.\tag 4.6
$$
First note that $H^0(X,\bbz_2) = \bbz_2$ for $X$ connected. One can verify (4.6) by choosing a section $\a$ of $U^\bot$ and computing that $\pi_*(\Div (\a )) \neq 0$ mod 2. (Note that  $\Div (\a )$ represents the mod 2 Euler class of $U^\bot$
which is equal to the top Stiefel--Whitney class, $w_{m-1} (U^\bot )$.)

An alternate proof of (4.6) can be given as follows. Equations (4.4) and (4.5)
imply that $\pi_* (w_{m-1} (U^\bot )) = \pi_* (a^{m-1} )$.
Using the standard fact that if $a$ is the nonzero element of $H^1 (\bbp (\bbr^m),\bbz_2)$ then $a^{m-1}$ is the non--zero element of $H^{m-1} (\bbp(\bbr^m),\bbz_2)$ we conclude that $\pi_*(a^{m-1} ) = 1$. 
\qed\enddemo

We only need Lemma 4.2 in the special case that $E = \ur^m$ is trivial.

\proclaim{Corollary 4.7}
Consider the tautological line bundle $U$ on $\bbp (\ur^m)$. Then
$$
\pi_*(w_1(U))^j \,\, = \,\, 0 \qquad\text{if } j \neq m-1 $$
and
$$\pi_* (w_1(U))^{m-1} = 1.$$
\endproclaim

\demo{Proof of Theorem 4.1}
As in Section 3, let $\nu : \ur^m\to F$ denote the bundle map corresponding to the sections $\nu_1,...,\nu_m$. Let $\bsnu : \ur^m\to\BF$ denote the pullback of $\nu$ to 
 the projectivization $\bbp (\ur^m)$. Let $\nut$ denote the restriction of 
$\bsnu$ to the tautological line bundle $U\subset\ur^m$. Then, considering $\Div ^{\oper{mod} 2}(\nut )\in \SF^{\oper{mod} 2}_\LOC (\bbp ( \ur^m))$ as a mod 2 current, the linear dependency current is defined to be the current pushforward
$$
\BL\BD ^{\oper{mod} 2}(\nu ) \,\, := \,\,\pi_* (\Div^{\oper{mod} 2} (\nut ))\in\SF^{\oper{mod} 2}_\LOC (X).
$$
Now, by Remark 2.17, the mod 2 divisor $\Div ^{\oper{mod} 2}(\nut )$ represents the top Stiefel-Whitney class of 
$H := \oper{Hom} (U,\BF )$ over $\bbp (\ur^m)$. The standard formula for the Stiefel--Whitney classes of a tensor product (see [MS]) says that
$$ w_n(H) \,\, = \,\, w_n (U^*\otimes \BF ) \,\, = \,\, \sum\limits_{j=0}^n w_{n-j} (\BF )\, (w_1 (U))^j.
$$
So, by Corollary 4.7, $\pi_*(w_n (H)) = w_{n-m+1} (F)$. Therefore $\BL\BD ^{\oper{mod} 2}(\nu )$ represents $w_q(F)$  as desired.
\qed\enddemo

Now we consider the  case that $q$ is odd  and identify the cohomology class of the $\SO_F$--twisted current $\BL\BD (\nu)\in\wt\SF_\LOC (X)$ in $H^q (X,\wt\bbz)$, where $\wt\bbz = \wt\bbz_F := \SO_F\otimes_{\bbz_2}\bbz$. 
Consider the short exact triple $0\to\wt\bbz\overset 2 \to \to \wt\bbz\to\bbz_2\to 0$ and the induced long exact sequence
$$
{\dots \to H^{q-1} (X,\bbz_2)\overset \b \to \to H^q (X,\wt\bbz) \overset 2 \to \to  H^q (X,\wt\bbz) \overset \rho \to \to  H^q (X,\bbz_2)\to \dots} \, .\tag 4.8
$$
Define the $q$th $\wt\bbz$--Stiefel--Whitney class $\wt W_q\in H^q (X,\wt\bbz)$
to be the Bockstein of $w_{q-1}$,
$$
\wt W_q := \b (w_{q-1}).\tag 4.9
$$
Recall from Example 7 of the Appendix that $\wt\bbz$--cohomology can be computed
using $\SO_F$--twisted integrally flat currents, $\wt\SF_\LOC (X)$.

\proclaim{Theorem 4.10. (q odd)}
Given an atomic collection $\nu$ of $m = n-q+1$ sections of a real vector bundle $F\to X$ the linear dependency current  $\BL\BD (\nu)$ represents the $q$th $\wt\bbz$--Stiefel--Whitney class, $\wt W_q(F)\in H^q (X,\wt\bbz)$, of the bundle $F$. Moreover, the mod 2 reduction of $\wt W_q$ equals $w_q$, i.e.  $\rho(\wt W_q) = w_q$. Hence ${\wt{\oper {Sq}}}^1 := \rho\circ\b$ maps $w_{q-1}$ to $w_q$ for $q$ odd.
\endproclaim

\demo{Proof}
A representative for $\wt W_q (F)$ can be computed as follows. Choose an atomic collection $\mu = (\mu_1,...,\mu_{m+1})$ of $m+1 = m- (q-1) +1$ sections. By Theorem 4.1, the mod 2 linear dependency current $\BL\BD^{\oper{mod} 2} (\mu)\in \SF_\LOC^{\oper{mod} 2} (X)$ represents the $(q-1)$th Stiefel--Whitney class 
$w_{q-1} (F)\in H^{q-1} (X,\bbz_2)$. Choose an $\SO_F$--twisted current representative $T\in\wt\SF_\LOC (X)$ of the mod 2 current  $\BL\BD^{\oper{mod} 2} (\mu)$. Then, by the definition of $\b$, the current $\tfrac 12 dT$ represents $\wt W_q (F)$.

Suppose that the subcollection $\eta =  (\mu_1,...,\mu_{m})$ is atomic and that
$\mu$ satisfies two additional assumptions described below. Then it is possible to choose $T$ so that 
$$
\tfrac 12 dT \,\, = \,\, \BL\BD (\eta)\qquad\text{on } X.\tag 4.11
$$
Consequently $ \BL\BD (\eta)$ represents $\wt W_q (F)$. The Theorem now follows from Theorem 3.7. Also note that, since the mod 2 reduction of $ \BL\BD (\eta)$
is $ \BL\BD^{\oper{mod} 2} (\eta)$, Theorem 4.1 implies that $\rho (\wt W_q) = w_q$.

The $\SO_F$-twisted current $T$ representing $\BL\BD^{\oper{mod 2}} (\mu)$ can be chosen as follows.
Firstly, if  $q=1$ define  $T$   to be the $\SO_F$--twisted generalized function $T := \tfrac {\s}{|\s |}$, where $\s$ is the section of $\SO_F^\bbr$ associated to 
$\eta$ (c.f. Remark 3.10). In this case (4.11) is simply a restatement of
Proposition 3.11.
Secondly, if  $q > 1$ is  odd,  $T$ is defined as follows. Let $\widehat \mu$ be the induced section of the rank $n$ bundle $\widehat H = \oper{Hom} (NS(\ur^{m+1}), \rho^*F)$ over $S(\ur^{m+1})$. Embed $\ur^m\hookrightarrow \ur^{m+1} = \ur^m\oplus\ur$ and let $t$ be the coordinate
on $\ur$. Set $\chi := \frac {t}{|t|}\in L^\infty_\LOC (S(\ur^{m+1}))$ and let
$p: S(\ur^{m+1})\to \bbp (\ur^{m+1})$ and $\pi: \bbp (\ur^{m+1})\to X$ denote the projection maps. Assume that the current $\Div (\widehat \mu )$ on $ S(\ur^{m+1})$ has locally finite mass. Then  $\Div (\widehat \mu )$ is an odd 
 locally rectifiable $\SO_F$--twisted current on $S(\ur^m)$  and so $\chi \Div (\widehat \mu )$ is a well defined even $\SO_F$--twisted flat current on $S(\ur^{m+1})$. Let  $\wt T := 
\tfrac 12 p_* (\chi \Div (\widehat \mu ))$ be the corresponding $\SO_{\bbp(\bbr^m)}$-twisted  current on $\bbp (\ur^{m+1})$. Since the mod 2 reduction of $\wt T$ 
 is $\Div^{\oper{mod} 2} (\wt\mu )$ its pushforward $T :=
\pi_* \wt T$ is an $\SO_F$--twisted current on $X$ whose mod 2 reduction is $\Div^{\oper{mod} 2} (\mu )$.

Next we verify that equation (4.11) holds for this choice of $T$. Suppose that
the codimension $n-1$ Hausdorff measure of $\oper{Zero}( \widehat \eta )$ is zero. Then, since $\widehat\eta = \widehat\mu\bigm|_{S(\ur^m)}$,
Lemma 4.12 below implies that
$$
d (\chi \Div (\widehat \mu )) \,\, = \,\, 2 i_* \oper{Div} (\widehat\eta )\qquad\text{on } S(\ur^{m+1}),
$$
where $i: S(\ur^m)\hookrightarrow S(\ur^{m+1})$ is the inclusion map. Finally,
by (3.6), this equation pushes forward to give equation (4.11) on $X$.
\qed\enddemo

\proclaim{Lemma 4.12}
Let $u: X\times \bbr\to \bbr^n$ be atomic and suppose that $\Div (u)$ has locally finite mass. Suppose that the function $v: X\to \bbr^n$ defined by $v(x) = u(x,0)$ is atomic and that the codimension $n-1$ Hausdorff measure of $\oper{Zero} (v)$ in $X$ is zero. Let $t$ denote the coordinate on $\bbr$ and define $i: X\hookrightarrow X\times \bbr$ by $i(x) = (x,0)$. Then
$$
d \( \tfrac{t}{|t|}\, \Div (u)\) \,\, = \,\, 2 i_* \Div (v) \qquad\text{on } X\times \bbr. \tag 4.13
$$
\endproclaim

\demo{Proof}
Let $\theta$ denote the solid angle kernel on $\bbr^n$. First note that, since
$i_* d = -\, di_*$, equation (4.13) is the exterior derivative of the degree $n$
current equation
$$
d\( \tfrac{t}{|t|} \, u^*\theta \) \,\, = \,\, \tfrac{t}{|t|}\, \Div (u) \,\, + \,\, 2 i_* (v^*\theta ) \qquad\text{on } X\times \bbr. 
$$
To verify this equation we argue as follows. Firstly, the equation holds on 
$(X\times\bbr)\sim(\oper{Zero} (v)\times \{ 0\})$ since it is true if $t\neq 0$, and if $t=0$ and $x\notin \oper{Zero}(v)$ then   $d(\tfrac {t}{|t|} \, u^*\theta )
= 2 [X] \, u^*\theta = 2 i_* (v^*\theta )$. Finally, since the codimension
$n$ Hausdorff measure of $\oper{Zero} (v)$ in $X\times\bbr$ is zero, 
the Federer Support Theorem for flat currents implies that the equation holds on all of $X\times\bbr$.
\qed\enddemo

\subheading{Remark 4.14. Non--injectivity currents}
Let $m=\oper{rk} E \leq \oper{rk} F =n$ and set $q=n-m+1$. In this remark we study the non--injectivity current of a bundle map $\nu: E\to F$. This is a 
degree $q$ current which is supported on the set of  points of $X$ over which the bundle map $\nu$ fails to be injective. 
 The  mod 2 non--injectivity current, $\BD_{\oper{NI}}^{\oper{mod} 2} (\nu)$, is defined by replacing $\ur^m$ by $E$  in Definition 3.14. The analogue of Remark 3.15 holds for  $\BD_{\oper{NI}}^{\oper{mod} 2} (\nu)$. Furthermore,  the cohomology class of  $\BD_{\oper{NI}}^{\oper{mod} 2} (\nu)$ in $H^q(X,\bbz_2)$ is $\{ w(F)\, w(E)^{-1} \}_q$, the degree $q$ part of $w(F)\, w(E)^{-1}$. 

 If $q$ is odd the   non--injectivity current, $\BD_{\oper{NI}}(\nu)$,
is defined as in Definition 3.4. (Note that, in the case that $q=n$ is odd, the current $\BD_{\oper{NI}}(\nu)$ is simply the divisor of the induced section of $\oper{Hom} (E,F)\to X$.) This current is an $\SO_E\otimes\SO_F$--twisted current on $X$. The analogues of all the results of Section 3A hold for $\BD_{\oper{NI}}(\nu)$. Furthermore the cohomology class of $\BD_{\oper{NI}}(\nu)$ in $H^q(X,\wt\bbz_{E\oplus F} )$ is $\b \( \{ w(F)\, w(E)^{-1} \}_{q-1} \)$, the Bockstein of the degree $q-1$ part of $w(F)\, w(E)^{-1}$.

\subheading{5. Obstruction currents}
\newline
\newline

The Stiefel--Whitney classes were  originally defined (see [S], [W]) as obstruction classes. 
The $q$th  {\bf obstruction class} of a real rank $n$ vector bundle $F\to X$
is a cohomology class which is the
 obstruction to the existence of a collection of $n-q+1$ linearly independent sections of $F$ over the $q$--skeleton of a cell decomposition of $X$. It is defined to be the cohomology class of a certain obstruction $q$--cocycle
which is associated to each suitable collection of $n-q+1$ sections of $F$. 
If this obstruction cocycle is defined for an atomic collection of sections, then, by Poincar\'e duality, there is also defined 
 a canonical {\bf obstruction current} on $X$. The aim of this section is to show that this 
  obstruction current is equal to the linear dependency current of this special  collection of sections.

We begin by recalling the definition of the {\bf obstruction cocycle} as given by Steenrod, [St]. Fix $q\in \{ 1,...,n\} $.  Let $\stF\to X$ be the bundle whose fibre over $x\in X$ is the Stiefel manifold consisting of all $(n-q+1)$--tuples of linearly independent vectors of $F_x$. 
Choose a smooth locally finite simplicial decomposition $K$ of $X$ and let $K'$
be the first barycentric subdivision of $K$. Each barycentrically subdivided $q$--simplex, $a^q$, of $K$ is a simplicial subcomplex of $K'$. In fact, since it is diffeomorphic to a $q$--ball, $a^q$ is a $q$--cell. If $q$ is odd or  $q=n$, choose an orientation on each cell.  The collection of such (oriented) cells forms a cellular subdivision $K_a$ of $K'$. Since $\pi_i (V_{n-q+1} (\bbr^n)) = 0$ for all $i < q-1$ there is a section $\eta$ of $\stF$ over the $(q-1)$--skeleton $K_a^{q-1}$ of $K_a$. 

Now  
$$\pi_{q-1} (V_{n-q+1} (\bbr^n))   \,\, = \,\,\, \left\{ 
\aligned 
\bbz\,\, & \qquad\text{if } q \text{ is odd or } q=n\\
\bbz_2 & \qquad\text{if } q \text{ is even and } q<n.
\endaligned
\right.
$$
Consequently there is no guarantee that $\eta$ will extend to a section of $\stF$ over the $q$--skeleton $K_a^q$. Steenrod defines a twisted (or mod 2) cellular $q$--cochain, $w_q(\eta )$, which is zero iff $\eta$ can be extended over $K_a^q$. Fix a point $x_a$ in each $q$--cell $a^q$ of $K_a^q$. The cochain $w_q(\eta )$ assigns an element $w_q(\eta ) (a^q)$ of $\pi_{q-1} (V_{n-q+1} (F_{x_a} ))$ to each $q$--cell $a^q$. It is defined as follows. Choose a trivialization of $F$ over $a^q$ and let $\psi: \stF\to  V_{n-q+1} (F_{x_a} )$ be the induced map. Then 
 $w^q(\eta ) (a^q)$ is defined to be the homotopy class of $\psi\circ\eta : \p a^q\to V_{n-q+1} (F_{x_a} )$. This class is well defined independent of the choice of trivialization of $F$. Steenrod shows that $w_q(\eta )$ is a cocycle
whose cohomology class in $H^q (X,\pist )$ is well defined independent of the choice of section $\eta$. By definition this class is the $q$th {\bf obstruction class} of $F$. In keeping with the notation of Section 4 the $q$th obstruction class will be denoted by $w_q(F)$ when $q$ is even and by $\wt W_q (F)$ when $q$ is odd.

In summary, if $q$ is odd or $q=n$ (resp. $q$ is even and $q<n$) Steenrod associates a $\pist$--twisted cochain (resp. mod 2 cochain), $w_q(\eta )$, to each section $\eta$ of the Stiefel bundle $\stF\to K_a^{q-1}$. On the other 
hand, in Section 3  we associated the  $\orF$--twisted current
$\BL\BD (\nu )$ (resp. mod 2 current $\BL\BD^{\oper{mod} 2} (\nu )$) to a collection of  $n-q+1$ sections of the vector bundle $F\to X$. Our
goal is to relate these two constructions. We begin by considering the case that $q$ is odd or $q=n$. 

\subheading{The case that $q$ is odd or $q=n$}
\newline
\newline

First note that there is a bundle isomorphism
$$\varphi : \pist\lra \wt\bbz_F$$
defined, in terms of a generator $[\s ]$ of $\pi_{q-1} (V_{n-q+1} (\bbr^n))$, as follows.
(See [St,25.6] for a definition of the homotopy generator $[\sigma]$.)
 Fix $x\in X$ and let $\eta : S^{q-1}\to V_{n-q+1} (F_x)$ represent
an element of $\pist_x$. For each frame $f$ of $F_x$ we obtain a map $\psi : V_{n-q+1} (F_x)\to V_{n-q+1} (\bbr^n)$. Define $\lambda\in \bbz$ by $[\psi\circ\eta ] = \lambda [\s ]$ in $\pi_{q-1}V_{n-q+1} (\bbr^n)$ and let $[f]$ denote the orientation class of the frame $f$ in $\SO_F$. Then $\varphi (\eta )$ is defined to be the class of $([f], \lambda)$ in $\wt\bbz_F$. (Recall that $\wt\bbz_F$ is the space of orbits of the $\bbz_2$--action $\rho ([f],\lambda )
= (-[f],-\lambda) $ on $ \SO_F\times \bbz$.)

Let $N = \oper{dim} X$. Now there
is a dual cellular decomposition $K_b$ of $K'$ characterized by the fact that to each $q$--cell $a^q$ of $K_a$ there is a unique $(N-q)$--cell, $b^{N-q}$, of $K_b$ so that the dual of each face of $a^q$ has $b^{N-q}$ as a face. Choose an orientation on each cell of $K_b$. Note that the intersection of a cell and its dual is the common centrepoint of both cells and that cells $a^q$ and $b^{N-q}$  which are not dual to each other do not intersect. Therefore the $(q-1)$--cells of $K_a$ do not intersect the $(N-q)$--cells of $K_b$. 

Using the fact that $\pi_i (V_{n-q+1} (\bbr^n)) = 0$ for $i < q-1$ we can construct  smooth sections $\nu_1,...,\nu_n$ of $V\to X$ so that for each $q\in \{ 1,...,n\}$ the linear dependency set of $\nu_1,...,\nu_{n-q+1}$ is a cellular subcomplex $K_b^{N-q} (\nu )$ of $K_b^{N-q}$. Note that $ \nu_1,...,\nu_{n-q+1}$
define a section $\eta$ of $\stF$ over $K_a^{q-1}$.

Now the cohomology group $H^*(X,\wt\bbz_F )$ can be computed using the smooth $\orF\otimes\orX$--twisted infinite integral dual cellular chain complex. 
The $(N-q)$--dimensional chains of this complex are formal infinite combinations of the form $\sum\limits_j \lambda_j\, b_j^{N-q}$ where the $\lambda_j$ are $\orF\otimes\orTX$--twisted integers and the $b_j^{N-q}$ are the oriented $(N-q)$--cells of $K_b$. Our intermediate goal is to associate to the cocycle $w_q (\eta )$ (defined by the collection
of sections $\nu_1, ..., \nu_{n-q+1}$) an $\orF\otimes\orTX$--twisted integral
cellular $(N-q)$--cycle, $w_{N-q} (\nu)$, which represents $\wt W_q(F)\in H^q(X;\wt\bbz_F )$.  It is defined as follows. First let
$a^q$ and $b^{N-q}$ be \od dual cells and set $\{ x \} := a^q\cap b^{N-q}$. Let $[a,b]$ denote the local section of $\orTX\to b^{N-q}$ induced by the decomposition $T_xX \cong T_x a^q \oplus T_x b^{N-q}$. Then
$$
w_{N-q} (\nu ) := \sum\limits_j \lambda_j\, b_j^{N-q}\tag 5.1
$$
 where the sum is taken over the cells $b_j^{N-q}$ of the linear dependency subcomplex $K_b^{N-q} (\nu )$ of $ \nu_1,...,\nu_{n-q+1}$, and where $\lambda_j$ is the $\orF\otimes\orTX$--twisted integer 
$$
\lambda_j := [a_j^q,b_j^{N-q}]\otimes \varphi( w^q(\eta)(a^q_j)).
$$
That is, $w_{N-q} (\nu )$ is the Poincar\'e dual of $w^q (\eta )$. 

Since an $\orX$-twisted smooth \od cellular chain defines a locally rectifiable current the 
twisted cellular cycle $w_{N-q} (\nu )$ defines an $\SO_F$--twisted
locally integrally flat current on $X$, which we also denote by $w_{N-q} (\nu )$. Now, by [HS,3.2], we can choose the sections $\nu_1,...,\nu_{n-q+1}$ so that the induced section $\widetilde\nu$ of $H\to\bbp (\ur^m)$ is atomic.
Then we have the following result.

\proclaim{Theorem 5.2}
For the collection of sections  $\nu_1,...,\nu_{n-q+1}$ described above,
$$
\BL\BD (\nu ) \,\, = \,\, w_{N-q} (\nu )
$$
as $\SO_F$--twisted
locally integrally flat currents on $X$. Consequently, the obstruction class $\wt W_q (F)\in H^q(X,\wt\bbz_F )$ is the
cohomology class of 
$\BL\BD (\nu)$. 
\endproclaim

\demo{Proof}
Let $m = n-q+1$. By construction $\nu_1,...,\nu_{m-1}$ are linearly independent
sections of $F$ over $X\sim K_b^{N-q-1}$. Let $E^{m-1}$ denote the \od span of $\nu_1,...,\nu_{m-1}$ over  $X\sim K_b^{N-q-1}$ and let $G^q$ denote the orthogonal complement of $E^{m-1}$ in $F$ with respect to some metric on $F$.
Let $\nu_{m}^\bot$ denote the orthogonal projection of $\nu_{m}$ onto $G^q$ over
$X\sim K_b^{N-q-1}$. Now we can choose   $\nu_1,...,\nu_{m-1}$ so that the section $\nu_m^\bot$ is atomic. Then we have the following
\enddemo

\proclaim{Lemma 5.3}
$$
 w_{N-q} (\nu ) \,\, = \,\, \Div (\nu_m^\bot )\qquad\text{over } X\sim K_b^{N-q-1}.
$$
\endproclaim

\demo{Proof of Lemma}
To prove the lemma we need to calculate each twisted integer $\lambda_j$ of equation (5.1) in terms
of the degree of a certain map between $(q-1)$--spheres. This was done by
Halperin and Toledo [HT] as follows.  Since $K_b^{N-q-1}$ and $K_a^q$ are disjoint the bundle $G^q$ is defined over the $q$--cell $a_j^q$. Choose a trivialization $\rho : G^q\to\bbr^q$ of $G^q$ over $a_j^q$. This induces an orientation on $G^q\to a_j^q$. Let $[f]$ denote the induced orientation on $F = E^{n-q}\oplus G^q$ over $a_j^q$. Then 
$$
\lambda_j = [a_j^q,b_j^{N-q}]\otimes [f] \,\, n_j
$$
 where $n_j$ is the degree of the induced map between \od $(q-1)$-spheres,
$$
n_j \,\, := \,\, \oper{Deg} \(\rho\circ \nu_{n-q+1}^\bot : \p a_j^q\to S(\bbr^q )\).
$$
The lemma now follows from Theorem  2.9 and Remark 2.11.
\qed\enddemo

\demo{Completion of the proof of Theorem 5.2}
By Theorem 3.18 and Lemma 5.3
$$
 w_{N-q} (\nu ) \,\, = \,\, \BL\BD (\nu )\qquad\text{over } X\sim K_b^{N-q-1}.
$$
Let $S :=  w_{N-q} (\nu ) - \BL\BD (\nu )$. Then $S$ is a flat current of dimension $N-q$ which is supported on the $N-q-1$ dimensional submanifold 
$K_b^{N-q-1}$. So, by the Federer Support Theorem for flat currents, [F], $S = 0$ on $X$, as 
required.
\qed\enddemo

\subheading{The case that $q$ is even and $q<n$}
\newline
\newline

In the case that $q$ is even and $q<n$ the obstruction class, $w_q(F)$, is an element of $H^q(X,\bbz_2)$. Now the cohomology group  $H^q(X,\bbz_2)$ 
can be computed using the smooth infinite mod 2 dual cellular chain complex. The
$(N-q)$--dimensional chains of this complex are formal infinite combinations of the form $\sum \lambda_j\, b_j^{N-q}$ where $\lambda_j\in\bbz_2$ and the $b_j^{N-q}$ are the (unoriented) $(N-q)$--cells of $K_b$. 

Just as in the case described above we construct sections $\nu_1,...\nu_{n-q+1}$
of $F$ and a mod 2 cellular $(N-q)$--cycle, $w^{\oper{mod} 2}_{N-q} (\nu ) = \sum \lambda_j\, b_j^{N-q}$, which is supported on the linear dependency set of  $\nu_1,...\nu_{n-q+1}$ and which represents the obstruction class $w_q(F)\in H^q(X,\bbz_2)$. This cycle defines a mod 2 locally integrally flat current
in the obvious way. Arguing as in the proof of Theorem 5.2 we have the following

\proclaim{Theorem 5.4}
For the collection of sections  $\nu_1,...,\nu_{n-q+1}$ described above,
$$
\BL\BD^{\oper{mod} 2} (\nu ) \,\, = \,\, w_{N-q}^{\oper{mod} 2} (\nu )
$$
as mod 2 locally integrally flat currents on $X$. Consequently, the obstruction class  $w_q(F)\in H^q(X,\bbz_2 )$ is the cohomology class of 
$\BL\BD^{\oper{mod} 2} (\nu)$.
\endproclaim

\subheading{Remark 5.5. Obstructions to injective bundle maps}
The results of this section can be generalized to the case of vector bundle
maps $\nu: E\to F$ where $m = \oper{rk} E \leq \oper{rk} F$. Let $\oper{Hom}^\times (E,F)\to X$ denote the bundle of injective bundle maps from $E$ to $F$. The Steenrod obstruction class of the bundle $\oper{Hom}^\times (E,F)\to X$ is a degree $q=n-m+1$ cohomology class which is the obstruction to the existence 
of an injective bundle map $\nu: E\to F$  over the $q$--skeleton
of $X$. Combining Remark 4.14 with the analogues of Theorems
5.2 and 5.4 we conclude that this class is equal to
\roster
\item the twisted Euler class $\wt e$ of $\oper{Hom} (E,F)$ when $q = n$
\item the degree $q$ part of $w(F)\, w(E)^{-1}$ when $q <n$ is even, and
\item the Bockstein of the degree $q-1$ part of $w(F)\, w(E)^{-1}$ when
$q<n$ is odd.
\endroster

\subheading{6. Higher dependency currents}
The aim of this section is to study the currents associated with higher dependencies, c.f. [HL3]. Let $\nu$ be an ordered collection of $m$ sections $\nu_1,...,\nu_m$ of a rank $n$ bundle $F\to X$. In Section 3 we studied the linear dependency current, $\BL\BD (\nu)$, which is supported on the set of points $x$ where at least one of the vectors $\nu_1 (x), ...,\nu_m (x)$ depend
linearly on the remaining ones. Fix an integer $\ell$ with $\oper{max} \{ 0, m-n\} \leq \ell \leq m$. In this section we study the higher dependency current, $\BL\BD_\ell (\nu)$,
which is supported on the set where at least $\ell$ of the sections depend 
linearly on the remaining ones, i.e. on the set of points over which the induced bundle map $\nu : \ur^m\to F$ has rank $\leq m-\ell$.
 The higher dependency current $\BL\BD_\ell (\nu)$ has degree $q := \ell ( n-m+ \ell )$, and so the dimension of $\BL\BD_\ell (\nu)$ decreases as
 $\ell$ increases.

The current  $\BL\BD_\ell (\nu)$ is defined as follows. Let $\pi : G_\ell (\ur^m)\to X$ denote the trivial Grassmann bundle of unoriented $\ell$--dimensional linear subspaces of the trivial bundle $\ur^m\to X$, and let $U\subset\ur^m$ be the tautological rank $\ell$ bundle over $G_\ell (\ur^m)$.
The collection $\nu$ is called
 $\ell$--{\bf dependency atomic}
 if the induced section $\wt\nu$ of the bundle $H = \oper{Hom} (U, \pi^*F)$ over $G_\ell (\ur^m)$ is atomic.

 The {\bf mod 2 higher dependency
current}, $\BL\BD_\ell^{\oper{mod} 2} (\nu)$, is defined by
$$
 \BL\BD_\ell^{\oper{mod} 2} (\nu) \,\, := \,\, \pi_* \Div^{\oper{mod} 2} (\wt\nu)\qquad\text{on } X.\tag 6.1
$$
In this mod 2 case orientation issues are irrelevant. 

If $m\equiv n$ mod 2 it is possible to pushforward the $\SO_H$-twisted current
 $\Div (\wt \nu )$ by the projection $\pi$. This is because
 there is a canonical isomorphism
$$
\SO_H \,\, \cong \,\, \cases \SO_{G_\ell (\bbr^m)}\otimes \pi^*\SO_F &\text{if $\ell$ is odd,} \\ \SO_{G_\ell (\bbr^m)} &\text{if $\ell$ is even.}\endcases
$$
 In this case the {\bf higher dependency current}, $\BL\BD_\ell (\nu)$, is defined by
$$
\BL\BD_\ell (\nu) \,\, := \,\, \pi_*\Div (\wt\nu ) \qquad\text{on } X.\tag 6.2
$$
Note that 
$$
\BL\BD_\ell (\nu ) \text{ is } \cases \text{ an } \SO_F-\text{twisted current on } X &\text {if $\ell$ is odd,} \\ \text{ a current on  $X$}  &\text{if $\ell$ is even.}\endcases
$$

The following equivalent definition of $\BL\BD_\ell (\nu)$ is often useful, (c.f. Remark 3.5). Let $\rho: \widehat G_\ell (\ur^m)\to X$ be the Grassmann
bundle of {\it oriented} $\ell$--dimensional linear subspaces of $\ur^m$. Note that the fibres $ \widehat G_\ell (\bbr^m)$ of $\widehat G_\ell (\ur^m)$ are 
canonically oriented manifolds. Let $\widehat U$ be the canonically oriented tautological bundle over $\widehat G_\ell (\ur^m)$, and let $\widehat \nu$ denote the induced section of the bundle $\widehat H = \oper{Hom} (\widehat U, \rho^*F)$ over $\widehat G_\ell (\ur^m)$. Then, if $m\equiv n$ mod 2,
$$
\BL\BD_\ell (\nu) \,\, = \,\, \tfrac 12 \, \rho_* \Div (\widehat \nu). \tag6.3
$$

\proclaim{Proposition 6.4}
Let $\nu :\ur^m\to F$ be as above and let $\psi :\ur^m\to\ur^m$ be a bundle isomorphism. Then the collection of
sections corresponding to the bundle map $\nu\circ\psi : \ur^m\to F$ is also
  $\ell$-- dependency atomic. Furthermore, if $n\equiv m$ mod 2,
$$
\BL\BD_\ell (\nu\circ\psi) \,\, = \,\, \cases \oper{sgn}\,\oper{det} (\psi)\,\,\BL\BD_\ell (\nu) &\text{if $\ell$ is odd,} \\ \BL\BD_\ell (\nu) &\text{if $\ell$ is even.}\endcases
$$
\endproclaim

\subheading{Note} The analogues of the main Theorems 3.7 and 3.15 for (mod 2)
linear dependency currents also hold for the (mod 2) higher dependency
currents. This fact together with Proposition 6.4 imply the following 
corollary.

\proclaim{Corollary 6.5}
Let $n\equiv m$ mod 2 and let  $\ell$ be {\bf odd}. Then there is a locally rectifiable $\SO_F$--twisted current $T$
so that
$$
2 \BL\BD_\ell (\nu) \,\, = \,\, dT.
$$
\endproclaim

\demo{Proof of  Proposition 6.4}
Let $\mu = \nu\circ\psi$ and let $\Psi : \widehat G_\ell (\bbr^m)\to \widehat G_\ell (\bbr^m)$ be the diffeomorphism induced by $\psi$. Arguing as in the proof of Proposition 3.8 it suffices to show that
$$
\Psi_* (\Div (\widehat\mu)) \,\, = \,\, \oper{sgn}\,\oper{det} (D\Psi) \, \Div(\widehat\nu) \qquad\text{on }  \widehat G_\ell (\bbr^m),\tag 6.6
$$
and that
$$
 \oper{sgn}\,\oper{det} (D\Psi) \,\, = \,\, \cases \oper{sgn}\,\oper{det} (\psi)  &\text{if $\ell$ is odd,} \\ 1 &\text{if $\ell$ is even.}\endcases \tag 6.7
$$
The proof of (6.6) is the same as that of Equation (3.13).
Next we prove (6.7). Clearly it suffices to consider orthogonal linear maps $\psi : \ur^m\to \ur^m$. Now there is a canonical orientation preserving bundle isomorphism $\phi: \oper{Hom} (\widehat U, {\widehat U}^\bot )\to T  \widehat G_\ell (\bbr^m)$ defined as follows. Fix $P\in  \widehat G_\ell (\bbr^m)$ and consider the canonical map $\varphi_P : \oper{Hom} (P, P^\bot) \hookrightarrow  \widehat G_\ell (\bbr^m)$ which sends a linear map to its graph. Let $\oper{Id} : \oper{Hom} (P, P^\bot) \to T_0 \oper{Hom} (P, P^\bot)$ be the canonical isomorphism. Then $\phi_P := D\varphi_P\circ \oper{Id}$ defines $\phi$ pointwise. Let $\overline{\Psi} : \oper{Hom} (P, P^\bot) \to \oper{Hom} (\psi(P), \psi(P^\bot))$ be the map defined by $\overline{\Psi} := \varphi_{\psi(P)}\circ\Psi\circ \varphi_P$. Now $\overline{\Psi} (\a) = \psi\circ\a\circ\psi^{-1}$ is a linear map and so $ \oper{sgn}\,\oper{det} (D\Psi)  = \oper{sgn}\,\oper{det} (\overline{\Psi})$. The result now follows by applying (3.3).
\qed\enddemo

\subheading{A. The mod 2 cohomology class of
 $\BL\BD_\ell^{\oper{mod 2}} (\nu )$}
The  goal of this subsection  is to identify the mod 2 cohomology class of the mod 2 higher dependency current. 
Let $w(F) = 1 + w_1(F) + w_2(F) + \dots + w_n (F)$ denote the total Stiefel--Whitney class of $F$. The {\bf Shur polynomial}, $\Delta^{(\ell)}_r (w(F))\in H^{rl} (X,\bbz_2)$, is the polynomial in $w_j(F)$ defined by
$$
\Delta^{(\ell)}_r (w(F)) \,\, := \,\, \oper{det} (w_{r-i+j} (F))_{1\leq i,j\leq\ell}. \tag 6.8
$$

\proclaim{Theorem 6.9}
Let $\nu:\ur^m\to F$ be as above and let
$q = \ell (n-m+\ell)$. Then the cohomology class of the mod 2 higher dependency current $\BL\BD_\ell^{\oper{mod} 2} (\nu)$ in $H^q(X,\bbz_2)$ is the Shur polynomial $\Delta^{(\ell)}_{n-m+\ell} (w(F))$.
\endproclaim

\subheading{Example 6.10}\newline
(a) Since $\Delta^{(3)}_1 (w(T\bbr\bbp^4)) \neq 0$, it is not possible to find a collection of 6 vector fields on $\bbr\bbp^4$ so that at each point of $\bbr\bbp^4$ at least 4 of the 6 vectors are linearly independent.
\newline
(b) Since $\Delta^{(3)}_3 (w(T\bbr\bbp^{10})) \neq 0$, for any collection of 10 vector fields on $\bbr\bbp^{10}$ there is a point of $\bbr\bbp^{10}$ so that at least 3 of the vectors at that point depend linearly on the remaining ones.

\demo{Proof of Theorem 6.9}
Let $\pi: G_\ell (\ur^m)\to X$. As in the proof of Theorem 4.1, it suffices to show that
$$
\pi_* w_{\ell n} (H) \,\, = \,\, \Delta^{(\ell)}_{n-m+\ell} (w(F))\qquad\text{in } H^q(X,\bbz_2).
$$
Now the standard formula for the Stiefel--Whitney class of a tensor product says that
$$
w_{\ell n} (H)\,\, = \,\, w_{\ell n} (U^*\otimes F) \,\, = \,\, \Delta^{(\ell)}_n
\( w(F)\, w(U)^{-1}\).
$$
(One way to see this is to apply equation (A.26) of [HL3] to Problem 7C of [MS].)
Let $k=m-\ell$. Then, arguing as in the proof of Theorem 4.4 of [HL3], we see that
$$
w_{\ell n} (H)\,\, = \,\,\Delta^{(\ell)}_{n-k} (w(F))\, \,\pi_* (w_k(\UP )^{\ell}).
$$
The proof is completed by observing that $w_k(\UP )^{\ell}$ is  the
generator of $H^{k\ell} (G_\ell (\bbr^{k+\ell}), \bbz_2) = \bbz_2$. 
\qed\enddemo

\subheading{Remark 6.11. Mod 2 degeneracy currents}
Let $\nu: E^m\to F^n$ be a bundle map. Fix an integer $k$ with $0\leq k < \oper{min} \{ m ,n\}$ and let $\ell = m-k$. The  mod 2 degeneracy current,
$\BD_k^{\oper{mod} 2} (\nu )$, of the bundle map $\nu$ is defined as in equation (6.1), with $\ur^m$ replaced by $E$. It is a degree $q = (m-k)(n-k)$ current which is supported on the set of points over which the bundle map has rank $\leq k$. The cohomology class of 
$\BD_k^{\oper{mod} 2} (\nu )$ in $H^q(X,\bbz_2) $ is given by $\Delta^{(m-k)}_{n-k} \( w(F)\, w(E)^{-1} \)$.
\newline
\newline

\subheading{Remark 6.12. Non--surjectivity currents}
Next we specialize Remark 6.11 to the case that $m= \oper{rk} E \geq \oper{rk} F
=n$ and  $k=n-1$ so that $q = m-n+1$. The  mod 2 non--surjectivity current, 
$\BD_{\oper{NS}}^{\oper{mod} 2} (\nu) : = \BD_{n-1}^{\oper{mod} 2} (\nu)$, is supported on the set over which the bundle map $\nu$ fails to be surjective. 
The cohomology class of  $\BD_{\oper{NS}}^{\oper{mod} 2} (\nu)$ in $H^q(X,\bbz_2)$ is given by $\Delta^{(q)}_{1} \( w(F)\, w(E)^{-1} \) \,\, = \,\, \Delta^{(1)}_{q} \( w(E)\, w(F)^{-1} \)\,\, = \,\, \{  w(E)\, w(F)^{-1}\}_q$, the degree $q$ part of $ w(E)\, w(F)^{-1}$. 

If $q$ is odd the  non--surjectivity current, $\BD_{\oper{NS}} (\nu) : = \BD_{n-1} (\nu)$, is an $\SO_E\otimes\SO_F$-twisted current
 defined as in equation (6.2).  Let $\nu^* : F^*\to E^*$ be the adjoint map. Then, at least for generic maps, $\BD_{\oper{NS}} (\nu)
= \BD_{\oper{NI}} (\nu^*)$, where $\BD_{\oper{NI}} (\nu^*)$ is the non-injectivity current of $\nu^*$ as defined in Remark 4.14. So, by Remark 4.14, the cohomology class of $\BD_{\oper{NS}} (\nu)$ in $H^q(X, \wt\bbz_{E\oplus F} )$ is $\b\(\{ w(E)\, w(F)^{-1}\}_{q-1}\)$, the Bockstein of the degree $q-1=m-n$ part of $w(E)\, w(F)^{-1}$.
\newline
\newline

\subheading{B. The integer cohomology class of $\BL\BD_\ell (\nu)$.
The case $n\equiv m$
 mod 2}
The aim of this subsection is to identify the (twisted) integer cohomology class $[\BL\BD_\ell (\nu)]$ of the higher dependency current $\BL\BD_\ell (\nu)$,
which is defined  whenever $n\equiv m$ mod 2.
The torsion free part of $[\BL\BD_\ell (\nu)]$ is well known. If $\ell$ is even
it  is a certain Shur polynomial in the total Pontrjagin class $p(F)$  (see
for example [HL3,6.9]) while if $\ell$ is odd it is zero, by Corollary 6.5
 above. The mod 2 reduction of $[\BL\BD_\ell (\nu)]$ is given 
in Theorem 6.9 above.  We will prove that the (twisted)
integer class of $\BL\BD_\ell (\nu)$ is the sum of its torsion free part and a
 2--torsion term, $T^{(\ell)}_{n-m+l} (\wt W)$, defined below. This result builds on work of Ronga [R] who showed that the integer class of $\BL\BD_{\ell}(\nu)$ is determined by its mod 2 and rational reductions. Our contribution is to explicitly identify the 2-torsion term as a certain polynomial in the Pontrjagin and 
twisted integral Stiefel-Whitney classes of $F$.

Throughout this subsection $\wt\bbz$ denotes the $\SO_F$-twisted integers, 
$\wt\bbz := \SO_F\otimes_{\bbz_2}\bbz$.
Let $p_i(F)\in H^i(X,\bbz)$denote the $i$th integral Pontrjagin class of $F$ and   $\wt W_{2j+1}(F)\in H^{2j+1} (X,\wt\bbz)$  the $(2j+1)$th twisted integral Stiefel--Whitney
class, defined by (4.9). Recall that $\rho (\wt W_{2j+1}) = w_{2j+1}$.  Note that, since $2\wt W_{2j+1} =0$, the subring of $H^*(X,\bbz)\oplus H^*(X,\wt\bbz)$ generated by the $p_i(F)$ and $\wt W_{2j+1} (F)$ is commutative. Also note that the product of two elements of $H^*(X,\wt\bbz)$ is an element
of $H^*(X,\bbz)$.
Let $\wt W_{2j}$ denote the
{\it formal} square root of the $j$th Pontrjagin class, i.e. $\wt W_{2j} := \sqrt{p_j}$. We make this definition because
$w^2_{2j}$ is the mod 2 reduction of $p_j$
  (see [MS]
Problem 15A). Of course the formal symbol $\wt W_{2j}$ has no cohomological meaning.

\subheading{Definition 6.13}
Let $S_\ell$ denote the symmetric group on $\ell$ elements. Define $\tau\in S_\ell$ by $\tau (i) = \ell +1-i$, and note that $\tau^2 = \oper{Id}$.  Let $R:S_\ell\to S_\ell$ be the involution defined by $R(\s) = \tau\s^{-1}\tau$. Define an index set $\SJ\subset S_\ell$ by
$$ 
\SJ \,\, := \,\,  \{\s\in S_\ell \,\, : \,\, R(\s) = \s \text{ and, if $\ell$ is even, then }  \s (i)\not\equiv i\oper{mod} 2 \text{ for some } i.\}
$$
Set $r:=n-m+\ell$ and note that $\ell\equiv r$ mod 2. Then we define
$$
T^{(\ell)}_r (\wt W) \,\,:= \,\, \sum\limits_{\s\in \SJ} \,\,\prod\limits_{i=1}^\ell
\wt W_{r+i-\s(i)}.
$$

\proclaim{Lemma 6.14}
Suppose $\ell\equiv r$ mod 2. Then\newline
(1) $T^{(\ell)}_r (\wt W)$ is a polynomial in the $p_i$  and $\wt{W}_{2j+1}$. 
\newline
(2) If $\ell$ is odd (resp. even) then each term of $T^{(\ell)}_r (\wt W)$ is of odd (resp. even) degree in the variables $\wt W_{2j+1}$. Therefore $T^{(\ell)}_r (\wt W(F))$ is an element of $H^q(X,\wt \bbz)$ (resp. $H^q(X,\bbz)$), where $q=\ell r$. 
\newline
(3) Each term of the polynomial $T^{(\ell)}_r (\wt W)$ has a factor of the form
$\wt W_{2j+1}$. Therefore the class $T^{(\ell)}_r (\wt W(F))$ is a torsion class of order 2.
\endproclaim

In the case that $\ell$ is even set $\ell = 2 \ell_0$ and $r=2r_0$.
 The main result is
 
\proclaim{Theorem 6.15}
The(twisted) integer cohomology class of the higher dependency current $\BL\BD_\ell (\nu)$ is 
$$
 [ \BL\BD_\ell (\nu)] \,\, = \,\, \cases  T^{(\ell)}_r (\wt W(F)) &\text{ in $H^q(X,\wt\bbz)$ when  $\ell$ is odd,}\\
\Delta^{(\ell_0)}_{r_0} (p(F)) \,\,+\,\, T^{(\ell)}_r (\wt W(F)) &\text{ in $H^q(X,\bbz)$ when  $\ell$ is even.}
\endcases
$$
\endproclaim

\subheading{Remark 6.16} For $\ell$ even (resp. $\ell$ odd) let  $\iota^* : H^*(X,\bbz)\to H^*(X,\bbr)$   (resp. $\iota^* : H^*(X,\wt\bbz)\to H^*(X,\wt\bbr)$) denote the usual coefficient
 homomorphism, and let 
$\rho^* : H^*(X, \bbz)\to H^*(X,\bbz_2) $ (resp. $\rho^* : H^*(X, \wt\bbz)\to H^*(X,\bbz_2) $) denote mod 2 reduction.
Then
$$
\iota^*  ( [ \BL\BD_\ell (\nu)]) \,\, = \,\, \cases  0 &\text{if $\ell$ is odd,}\\
\Delta^{(\ell_0)}_{r_0} (p(F)) &\text{if $\ell$ is even.}
\endcases\tag 6.17$$
and 
$$
\rho^* ( [ \BL\BD_\ell (\nu)]) \,\, = \,\, 
 [ \BL\BD_\ell^{\oper{mod} 2} (\nu)] 
\,\, = \,\, \Delta^{(\ell)}_r (w(F)). \tag 6.18
$$
These two formulae are well known (see [R]) and will be used to prove
 the more general result of Theorem 6.15.

\subheading{Example 6.19}
\newline
(1) If $\ell =2$ then $r = 2r_0 = n-m+2$ and $[\BL\BD_2(\nu)] = p_{r_0} + \wt W_{r-1} \wt W_{r+1}$.
\newline
(2) If $\ell = r = 3$ then
$
[\BL\BD_3 (\nu)] \,\, = \,\, p_1\wt W_5 \,\, + \,\, p_2\wt W_1\,\, + \,\, \wt W_3^3 \,\, + \,\, \wt W_1 \wt W_3\wt W_5.
$\newline
(3) If $\ell = r = 4$  then
$[\BL\BD_4(\nu)]$ 
$$ =  p_2^2-p_1p_3+p_1\wt W_5\wt W_7+p_2(\wt W_3\wt W_5+\wt W_1\wt W_7) +p_3\wt W_1\wt W_3+(\wt W_1\wt W_5+\wt W_3^2)(\wt W_3\wt W_7+\wt W_5^2).
$$

\subheading{Remark 6.20} The following equivalent definition of $T^{(\ell)}_r (\wt W)$ will be useful. 
Define an $\ell\times\ell$ matrix $(a_{ij})$ by
$$
a_{ij} \,\, := \,\, \wt W_{r+i-j} \qquad 1 \leq i,j\leq \ell.  \tag 6.21
$$
Note that the matrix $(a_{ij})$ is symmetric under reflection in the antidiagonal $i+j=\ell +1$, i.e. $a_{ \ell +1 -j\, ,\, \ell +1 -i} = a_{ij}$.
Let $\SI$ denote the collection of subsets $I$  of the index
set $\{ (i,j) \,\, : \,\, 1\leq i,j\leq\ell \}$ which satisfy the following three
properties:
\roster
\item
For each $i$ (resp. $j$) in $\{ 1,...,\ell\}$ there there is exactly one element
$j$ (resp. $i$) of $\{ 1,...,\ell\}$ so that $(i,j)\in I$.
\item The subset $I$ is symmetric under reflection in the antidiagonal $i+j=\ell +1$, i.e. $(i,j)\in I$ iff $(\ell +1 -i, \ell +1-j)\in I$, and
\item
If $\ell$ is even, then there is at least one element $(i,j)\in I$ for which
$i\not\equiv j$ mod 2.
\endroster

Then 
$$
T^{(\ell)}_r (\wt W) \,\, = \,\, \sum\limits_{I\in \SI} \,\, \prod\limits_{(i,j)\in I} \wt W_{r+i-j}. \tag 6.22
$$
To see that Definition 6.13 and (6.22) agree note that  each $\s\in\SJ$
defines a subset $I(\s)$ of $\{ (i,j) \,\, : \,\, 1\leq i,j\leq\ell\} $ by
$I(\s) := \{ (i,\s(i) ) \,\, : \,\, 1\leq i\leq \ell\} $. Furthermore the set associated with $R(\s)$ is the reflection in the line $i+j=\ell +1$ of the subset associated with $\s$.

\demo{Proof of Lemma 6.14}
 Let $\s\in\SJ$. Then, since $R(\s) = \s$, we have $(i,j)\in I(\s)$ iff $(\tau(j),\tau(i))\in I(\s)$, where $\tau (i) = \ell +1-i$. So, since $a_{\tau(j) \tau(i)} = a_{ij}$, 
$$
T^{(\ell)}_r (\wt W) \,\, = \,\, \sum\limits_{\s\in \SJ} \prod\limits_{i\in U(\s)} \wt W^2_{r+i-\s(i)} \prod\limits_{i\in\Delta (\s) }\wt W_{r+i-\s(i)},
$$
where $U(\s) = \{ i\,\, : \,\, i+\s(i) < \ell +1\}$ and $\Delta (\s) = \{ i\,\, : \,\, i+\s(i) = \ell +1\}$. Note that $|\Delta (\s)| \equiv \ell$ mod 2. Conclusions (1,2,3) now 
follow from the fact that, if $i\in \Delta(\s)$, then $r+i-\s(i) = r-\ell +2i-1$ is odd. Conclusion (4) follows from the fact that $|\Delta (\s)|\neq 0$ in the case that $\ell$ is odd, and from the definition of $\SJ$ in the case that $\ell$ is even. 
\qed\enddemo

\demo{Proof of Theorem 6.15}
By naturality we can reduce to the case in which the bundle $F\to X$ is the tautological rank $n$ bundle $U$ over a sufficiently high dimensional approximation, $G_n (\bbr^N)$, to the classifying space $G_n(\bbr^\infty)$.
Now, if $N$ is large enough, the torsion subgroup of $H^q(G_n(\bbr^N),\bbz)$ is a direct sum of cyclic groups of order 2 (see [B]). Furthermore, choosing 
 $N$ to be  odd,
$\SO_{TG_n(\bbr^N)} \cong \SO_U$ and so $H^*(G_n(\bbr^N), \wt\bbz) \cong H_*(G_n(\bbr^N),\bbz)$. So, by the universal coefficient theorem, the torsion subgroup of $H^q(G_n(\bbr^N),\wt\bbz)$ is also a direct sum of cyclic groups of order 2.
Consequently elements of $H^q(G_n(\bbr^N),\bbz)$ and $H^q(G_n(\bbr^N),\wt\bbz)$
are completely determined by their mod 2 and real reductions. So, setting
$$
Q^{(\ell)}_r \,\, := \,\, \cases  T^{(\ell)}_r (\wt{W}(F)) &\text{if $\ell$ is odd,}\\
\Delta^{(\ell_0)}_{r_0} (p(F)) \,\,+\,\, T^{(\ell)}_r (\wt{W}(F)) &\text{if $\ell$ is even,}
\endcases
$$
it suffices to prove that
$$
\iota^* (Q^{(\ell)}_r) \,\, = \,\, \cases  0 &\text{if $\ell$ is odd,}\\
\Delta^{(\ell_0)}_{r_0} (p(F)) &\text{if $\ell$ is even,}
\endcases\tag 6.23
$$
and 
$$
\rho^*(Q^{(\ell)}_r) \,\, = \,\, \Delta^{(\ell)}_r (w(F)),\tag 6.24
$$
that is, that $Q^{(\ell)}_r$ and $[\BL\BD_\ell (\nu)]$ have the same 
torsion free part and 
mod 2  reduction.

Now (6.23) follows immediately from Lemma 6.14 (3). 
Let $w(\s) := \prod\limits_{i=1}^\ell w_{r+i-\s(i)}$. Then, since the 
matrix $a_{ij} = w_{r+i-j}$ is symmetric under reflection in the antidiagonal, 
$$
w(R(\s)) \,\, = \,\, w(\s). \tag 6.25
$$
In the case that $\ell$ is odd we verify (6.24) by observing that, by Definition 6.13,
$
\Delta^{(\ell)}_r (w(F))\,\, - \,\, \rho^*(Q^{(\ell)}_r) \,\, = \,\, \sum\limits_{
\s\notin  \SJ} w(\s),
$
which is zero, since, by (6.25),  the sum is a sum of terms of the form $w(\s) + w(R(\s)) = 2 w(\s) = 0$.

Finally we verify (6.24) in the case that $\ell = 2\ell_0$ is even. Let $\psi : S_{\ell_0}
\to S_\ell$ be the injection defined for $j\in \{ 1,2,...,\ell_0\} $ by
$$
\psi(\eta) (2j-1) := 2\eta (j)-1\qquad\text{and}\qquad \psi(\eta) (2j) := 2\eta (j).
$$
The map $\psi$ can be interpreted as follows. Let $(b_{ij})$ denote the $\ell_0\times\ell_0$ matrix $b_{ij} := p_{r_0 +i-j}$ and let $C(\eta) =\{ b_{j\, ,\,\eta(j)} \,\, : \,\, 1\leq j\leq \ell_0\}$ be the subset of entries of $(b_{ij})$ defined
by $\eta\in S_{\ell_0}$. Then $C(\psi(\eta)) =\{ a_{i\, ,\,\psi(\eta)(i)} \,\, : \,\, 1\leq i\leq \ell \}$ is the set of those entries of $(a_{ij})$ obtained from $C(\eta)$
by replacing each element $b_{j\, ,\,\eta(j)} $ of $C(\eta)$ by the diagonal entries of the corresponding $2\times 2$ submatrix
$\left(\smallmatrix w_{2(r_0 +j -\eta(j))} & 0 \\
        0 & w_{2(r_0 +j -\eta(j))}\endsmallmatrix\right)$
of the matrix $(a_{ij})$. So, since $\rho^*(p_j) = w^2_{2j}$, the mod 2 reduction of the product of the elements of $C(\eta)$ equals the product, $w(\psi(\eta))$,
of the elements of $C(\psi(\eta))$. Summing over $\eta\in S_{\ell_0}$ we conclude that
$$
\rho^*\(\Delta^{(\ell_0)}_{r_0} (p(F))\) \,\, = \,\, \sum\limits_{\s\in \psi(S_{\ell_0})} w(\s).
$$

Now let $\SK := S_\ell\sim (\psi(S_{\ell_0}) \cup \SJ)$. To verify (6.24) it suffices to show that
$$
\sum\limits_{\s\in\SK} w(\s) \,\, = \,\, 0.\tag 6.26
$$
To prove (6.26) we study the index set $\SK$. Define $\b\in S_\ell$ by
$$
\b (2j-1) \,\, = \,\, 2j\qquad\text{and}\qquad \b (2j) \,\, = \,\, 2j-1\qquad\text{for } j\in \{ 1,...,\ell_0\},
$$
and define $P: S_\ell\to S_\ell$ by $P(\s) := \b\s\b$. 
The involution $P$ can be interpreted as follows. Firstly, each entry $a_{ij}$ of the $\ell\times\ell$ matrix $(a_{ij})$ has a pair $P(a_{ij})$ defined as follows. Partition $(a_{ij})$ into  $2\times 2$ submatrices. Let
$\left(\smallmatrix a&b\\c&d\endsmallmatrix\right)$ be one such submatrix. Then $P(a) = d$ and $P(b) =c$. Let $C(\s)$ denote the set of entries of $(a_{ij})$ defined by $\s\in S_\ell$. Then $C(P(\s)) = C(\s)$.
The pairing involution $P$ is introduced because
$$
\psi (S_{\ell_0}) \,\, =  \,\, \{ \s\in S_\ell \,\, : \,\, P(\s) =\s \text{ and } i\equiv \s(i) \oper{mod} 2 \text{ for all } i\}.\tag 6.27
$$
Define
$$
\SK_R = \{ \s\in \SK \,\, : \,\, R(\s) \neq\s\}\qquad\text{and}\qquad
\SK_P = \{ \s\in \SK \,\, : \,\, R(\s) =\s \text{ and } P(\s)\neq \s\}.
$$
We claim that $\SK$ is the disjoint union $\SK = \SK_R\cup\SK_P$. To see this choose $\s\in \SK \sim\SK_R$. Then, since $\s\not\in \SJ$, $i\equiv\s(i) $ mod 2 for all $i$. Therefore, since $\s\not\in \psi (S_{\ell_0})$, equation (6.27)
implies that $P(\s) \neq \s$, as required.

Note that, since the pair of the reflection of an entry of $(a_{ij})$ is the reflection of the pair of that entry,
$$
R(P(\s)) \,\, = \,\, P(R(\s)).\tag 6.28
$$
Hence $R$ preserves the decomposition $S_\ell = \psi(S_{\ell_0})\cup \SJ\cup\SK_R\cup\SK_P$. Now, since the involution $R:\SK_R\to\SK_R$ has no fixed points, it follows that $\sum\limits_{\s\in\SK_R} w(\s) =0$, since it is the sum
of terms of the form $w(\s) + w(R(\s)) = 2 w(\s) =0$. Finally, by (6.28), 
$P:\SK_P\to\SK_P$ is an involution with no fixed points, and once again  $\sum\limits_{\s\in\SK_P} w(\s) =0$. Hence (6.26) holds, as desired.
\qed\enddemo

\subheading{7. Applications}

In this section we apply the general results of the previous sections to study singularities of projections and singularities of maps (c.f. [HL3]).
 The results of this section hold whenever the projections and maps in question are atomic,
by which we mean that the induced section of $\oper{Hom} (U,F)\to G_\ell (E)$ is atomic.  This hypothesis is assumed throughout. In particular, in the real analytic case we simply require that the degeneracy subvarieties of the map have codimension greater than or equal to the
expected codimension in $X$, (see [HL3, 2.14]).

\subheading{A. Singularities of Projections}

Let $j: X\to \bbr^N$ be an immersion of a smooth $m$--manifold. Fix an integer $n < N$ and let $P: \bbr^N\to\bbr^n$ be a linear map. We study the singularities of the smooth projection $\widehat P = P\circ j : X\to \bbr^n$. Fix an integer
$k$ with $0\leq k < \oper{min} \{ m,n\} $. The $k$th {\bf mod 2 degeneracy current} of the projection $P$ on $X$ is defined to be $\BD_k^{\oper{mod} 2} (P) := \BD_k^{\oper{mod} 2} (d\widehat P)$ (c.f. Remark 6.11). This current is a 
degree $q:= (m-k)(n-k)$ current which is supported on the set where the 
differential $d\widehat P : TX\to\bbr^n$ has rank $\leq k$. By Remark 6.11,
$$
[\BD_k^{\oper{mod} 2} (P)] \,\, = \,\, \Delta^{(m-k)}_{n-k} (w(TX)^{-1}) \,\, = \,\, \Delta^{(n-k)}_{m-k} (w(TX))\qquad\text{in } H^q(X,\bbz_2). \tag 7.1
$$

If $n\equiv m$ mod 2, the $k$th {\bf degeneracy current} of the projection $P$ on $X$, $\BD_k (P) := \BD_k (d\widehat P)$ can also be defined (as in (6.2)). Let $(d{\widehat P})^* : \bbr^n \to T^*X$ denote the adjoint map. Then, at least for generic $P$, $\BD_k (P) = \BD_k ((d{\widehat P})^*)$. Therefore,
by Theorem 6.15,
$$
 [ \BD_k (P)] \,\, = \,\, \cases  T^{(n-k)}_{m-k} (\wt W(TX)) &\text{ in $H^q(X,\wt\bbz)$, when  $n-k$ is odd,}\\
\Delta^{(n_0 - k_0)}_{m_0-k_0} (p(TX)) \,\,+\,\, T^{(n-k)}_{m-k} (\wt W(TX)) &\text{ in $H^q(X,\bbz)$, when  $n-k$ is even,}
\endcases
$$
where $2(n_0-k_0) = n-k$ and $2(m_0-k_0) = m-k$.

\subheading{Example 7.3. Tangential Stiefel--Whitney classes}
Fix $1 \leq q\leq m$ and let $P: \bbr^N\to\bbr^{m-q+1} $ be  linear. The mod 2 
{\bf non--submersion current} of the projection $P$ on $X$ is defined by 
$\BD_{\oper{NS}}^{\oper{mod} 2} (P) := \BD_{m-q}^{\oper{mod} 2} (d\widehat P)$ on $X$. This is a degree $q$ current which is supported on the subset of $X$
on which the map $\widehat P : X\to \bbr^{m-q+1}$ fails to be a submersion.
By (7.1),
$$
[\BD_{\oper{NS}}^{\oper{mod} 2} (P)] \,\, = \,\, w_q(TX)\qquad\text{in } H^q(X,\bbz_2).
$$
Furthermore, if $q$ is odd, then the  non--submersion current, $\BD_{\oper{NS}} (P)$, can also be defined, and, by Remark 6.12,
$$
[ \BD_{\oper{NS}} (P)] \,\, = \,\, \wt W_q (TX)\qquad\text{in } H^q(X,\wt\bbz).
$$

\subheading{Example 7.4. Normal Stiefel--Whitney classes}
Fix $1 \leq q\leq N-m$ and let $P: \bbr^N\to\bbr^{m+q-1} $ be  linear. The mod 2 {\bf non--immersion current} of the projection $P$ on $X$ is defined by 
$\BD_{\oper{NI}}^{\oper{mod} 2} (P) := \BD_{m-1}^{\oper{mod} 2} (d\widehat P)$ on $X$. This is a degree $q$ current which is supported on the subset of $X$
on which the map $\widehat P : X\to \bbr^{m+q-1}$ fails to be a immersion.
By (7.1),
$$
[\BD_{\oper{NI}}^{\oper{mod} 2} (P)] \,\, = \,\, w_q(NX)\qquad\text{in } H^q(X,\bbz_2),
$$
where $NX$ is the normal bundle to $X$ in $\bbr^N$.
Furthermore, if $q$ is odd, then the  non--immersion current, $\BD_{\oper{NI}} (P)$, can also be defined, and, by Remark 4.14,
$$
[ \BD_{\oper{NI}} (P)] \,\, = \,\, \wt W_q (NX)\qquad\text{in } H^q(X,\wt\bbz).
$$

\subheading{B. Singularities of maps}
\newline
\newline

Let $X$ and $Y$ be smooth manifolds of dimensions $m$ and $n$ respectively, and
let $f:X\to Y$ be a smooth mapping. Let $0\leq k < \oper{min} \{ m,n\} $. The
$k$th mod 2 degeneracy current of the map $f$ is defined to be $\BD_k^{\oper{mod} 2} (f) := \BD_k^{\oper{mod} 2} (df)$. This is a degree
$q=(m-k)(n-k)$ current supported on the set where $df: TX\to TY$ has rank $\leq k$. By Remark 6.11, 
$$
[\BD_k^{\oper{mod} 2} (f)] \,\, = \,\, \Delta^{(m-k)}_{n-k} \( f^*(w(TY))\, w(TX)^{-1}\)\qquad\text{in } H^q (X,\bbz_2). \tag 7.5
$$

\subheading{Example 7.6. Non--submersion currents}
Suppose that $m=\oper{dim} X \geq \oper{dim} Y =n$ and let $q = m-n+1$. Then, by Remark 6.12, the cohomology class of the mod 2 non--submersion current $\BD_{\oper{NS}}^{\oper{mod} 2} (f) := \BD_{n-1} ^{\oper{mod} 2} (f)$ is
$$
[\BD_{\oper{NS}}^{\oper{mod} 2} (f)] \,\, = \,\, \{ w(TX)\, f^*(w(TY)^{-1})\}_{\oper{deg} q} \qquad\text{in } H^q(X,\bbz_2),
$$
and, if $q$ is odd,
$$
[\BD_{\oper{NS}} (f)] \,\, = \,\, \b\(\{ w(TX)\, f^*(w(TY)^{-1})\}_{\oper{deg} q-1}\) \qquad\text{in } H^q(X,\wt\bbz).
$$

\subheading{Example 7.7. Non--immersion currents}
Suppose that $m=\oper{dim} X \leq \oper{dim} Y =n$ and let $q = n-m+1$. Then, by Remark 4.14, the cohomology class of the mod 2 non--immersion current $\BD_{\oper{NI}}^{\oper{mod} 2} (f) := \BD_{m-1} ^{\oper{mod} 2} (f)$ is
$$
[\BD_{\oper{NI}}^{\oper{mod} 2} (f)] \,\, = \,\, \{ f^*w(TY)\, w(TX)^{-1}\}_{\oper{deg} q} \qquad\text{in } H^q(X,\bbz_2),
$$
and, if $q$ is odd,
$$
[\BD_{\oper{NI}} (f)] \,\, = \,\, \b\(\{ f^*w(TY)\, w(TX)^{-1}\}_{\oper{deg} q-1}\) \qquad\text{in } H^q(X,\wt\bbz).
$$

\subheading{Appendix. Computing cohomology with currents}

This appendix is included for two reasons; first for the sake of completeness.
The second reason is that although the  approach taken here is both simple and natural (via standard sheaf theory)  it does not appear in the geometric measure theory literature. 

\subheading{Definition A.1}
The complex 
$$0\to\SS\to\SF^0\overset d \to \to \SF^1\overset d \to \to \dots 
\overset d \to \to \SF^n\to 0$$
of sheaves is called an {\bf acyclic resolution of the sheaf} $\SS$ if
\roster
\item the complex is exact, and
\item each sheaf $\SF^p$ is acyclic, i.e. $H^j(X, \SF^p ) = 0 $ for $ j = 1,...\,\,.$ 
\endroster
\newline
\newline

The basic result is that cohomology with coefficients in $\SS$ can be computed from such a resolution. That is,
$$
H^p(X,\SS ) \,\, = \,\, \f{\{\varphi\in\Gamma (X,\SF^p ) \,\, : \,\, d\varphi =0\}}{d\Gamma (X, \SF^{p-1})}
$$
and 
$$
H^p_{\oper{cpt}}(X,\SS ) \,\, = \,\, \f{\{\varphi\in\Gamma_{\oper{cpt}} (X,\SF^p ) \,\, : \,\, d\varphi =0\}}{d\Gamma_{\oper{cpt}} (X, \SF^{p-1})}.
$$
A classical reference for this and other standard results from sheaf theory is Godement [G]. 

In this paper the cases and coefficient sheaves of most interest are:
\newline
{\bf The integer case} with coefficient sheaf $\bbz$, the sheaf of germs of locally constant integer
valued functions (see Corollary A.5).
\newline
{\bf The mod 2 case} with coefficient sheaf $\bbz_2 := \bbz \bigm / 2\bbz $
(see Example 10), and
\newline
{\bf The twisted integer case} with coefficient sheaf $\wt\bbz := \bbz\otimes_{\bbz_2} \SO_V$, where
$\SO_V$ is the orientation sheaf of a real vector bundle $V\to X$ (see Example 7).

Throughout this appendix $X$ is a $C^\infty$ (paracompact) $n$--dimensional manifold. Let $\SO_{TX}$ or $\SO_X$ denote the orientation sheaf of $X$ and let $\wt\bbz_X := \bbz\otimes_{\bbz_2} \SO_X$ and $\wt\bbr_X := \bbr\otimes_{\bbz_2}\SO_X$.

\subheading{Example 1. Differential forms}
Let $\SE$ denote the sheaf of germs of $C^\infty$ differential $p$--forms with $d$ taken to be exterior differentiation. Then 
$$
H^p(X,\bbr ) \,\, = \,\, \f{\{\varphi\in\SE^p(X)  \,\, : \,\, d\varphi =0\}}{d\SE^{p-1} (X)}.
$$
Each $\SE^p$ is acyclic because $\SE^p$ is fine (i.e. there exists a partition 
of unity). Also the sequence,
$
0\to\bbr\to\SE^0\to\dots\to\SE^n\to 0
$
is exact by the Poincar\'e lemma for exterior differentiation.
\newline
\newline

\subheading{Example 2. Currents}
Let ${{\SD'}}^p$ denote the sheaf of germs of degree $p$ currents (defined in
Section 2). Then
$$
0\to\bbr\to{{\SD'}}^0\to{{\SD'}}^1\to \dots \to
{{\SD'}}^n\to 0
$$
is an acyclic resolution of $\bbr$ and hence can be used to compute real cohomology $H^p(X,\bbr )$.
\newline
\newline

\subheading{Example 3. Singular chains}
Let $\SC_k$ denote the sheaf of germs
of $C^\infty$ singular $k$--chains with locally finite support, with the usual boundary operator. Then
$$
0\to\wt\bbz_X\to \SC_n\to\SC_{n-1}\to\dots\to\SC_0\to 0
$$
is an acyclic resolution of $\wt\bbz_X$, so that 
$$
H^p(X,\wt\bbz_X )\,\, = \,\, \f{\{\varphi\in{\SC}_k(X)  \,\, : \,\, d\varphi =0\}}{d{\SC}_{k+1} (X)}, \qquad\text{where } p + k =n.
$$
Homology is just cohomology with compact supports and with the coefficients 
twisted by the orientation sheaf $\SO_X$. That is,
$$
H_k(X,\bbz ) \,\, = \,\, H^p_{\oper{cpt}}(X,\wt\bbz_X )\,\, = \,\, \f{\{\varphi\in{\SC}^{\oper{cpt}}_k(X)  \,\, : \,\, d\varphi =0\}}{d{\SC}^{\oper{cpt}}_{k+1} (X)}, \qquad\text{where } p+k =n.
$$
\newline

\subheading{Example 4. $\orX$-twisted currents}
Let $\SD'_k$  denote the sheaf of germs of $k$ dimensional $\orX$-twisted
 currents on $X$. Then
 $0\to\wt\bbr_X\to \SD'_n\to\SD'_{n-1}\to\dots\to\SD'_0\to 0$ is an acyclic
resolution of $\wt\bbr_X$
 and hence can be used to 
compute $H^p(X, \wt\bbr_X )$ or real homology $H_k(X,\bbr ) := H^p_{\oper{cpt}} (X,\wt\bbr_X )$, where $p+k=n$.  
\newline
\newline

\subheading{Example 5. Integrally flat currents}
This is one of the  examples of central importance in this paper and so will be treated in more detail.  Let $\SF^p_\LOC (X)$ 
denote the space of locally integrally flat degree $p$ currents
on $X$. 
 We take as definition $\SF^p_\LOC (X) := \SR^p_\LOC (X) + d\SR^{p-1}_\LOC (X)$
i.e. all currents which can be written as $A + dB$ with $A\in \SR^p_\LOC (X)$ and $B\in \SR^{p-1}_\LOC (X)$ where $\SR^p_\LOC (X)$ denotes the space of locally rectifiable degree $p$ currents.
 In results where the  degree of a locally integrally flat current can be arbitrary we use the less encumbered notation $\SF_\LOC (U)$.

The spaces $\{\SF_\LOC (U)\,\, : \,\, U^{\oper{open}} \subset X\}$ form
a presheaf of abelian groups. One can form the associated sheaf $\SF_\LOC$ of germs, and consider the natural map from $\SF_\LOC (U)$ to $\Gamma (U, \SF_\LOC )$. This map is injective because the support axiom is satisfied.

\subheading{Support axiom} Let $\SF$ be a presheaf. If $T\in \SF (X)$ restricts to be zero in a neighbourhood of each point of $X$ then $T =0$.

\subheading{Note} For a general presheaf $\SF$ this axiom is equivalent to the 
concept of support being well defined. The support of $T\in \SF (X)$ is defined to be the complement of the set of points $x\in X$ such that $T\bigm|_U = 0$ for
some neighbourhood $U$ of $x$. 
\newline
\newline
This map is surjective because the following axiom is satisfied.

\subheading{Local to global axiom} 
Suppose $T_\a\in\SF(U_\a )$ is given, where $\{ U_\a\}$ is a locally finite open cover of $X$. Let $U_{\a\b} = U_\a\cap U_\b$. If $T_{\a\b} = T_\a\bigm|_{U_{\a\b}} - \,\, T_\b\bigm|_{U_{\a\b}}$ vanishes then there exists a global $T\in\SF (X)$ such that $T\bigm|_{U_\a} = \,T_\a$. 
\newline
\newline

For a given presheaf, if both of these conditions/axioms are satisfied then the presheaf is said to be a {\bf sheaf}. A sheaf is said to be {\bf soft} if for 
each closed set $C\subset X$ and each section of the sheaf on $C$ there exists an extension to all of $X$. That is, for each section on a neighbourhood of $C$ there exists a section on $X$ which agrees with the given section on a (smaller)
neighbourhood of $C$. Soft sheaves are always acyclic, for basically the same reason fine sheaves are acyclic. Namely, the decompositions provided by a partition of unity exist (even though these may  not arise from a partition of unity).

\proclaim{Theorem A.2}
The presheaf $\{\SF_\LOC (U)\}$  of locally integrally flat currents
 is a sheaf and this sheaf $\SF_\LOC$ is soft.
\endproclaim

\demo{Proof}
Since each $\SF_\LOC (U)$ is a subset of ${{\SD'}}(U)$ and the support axiom is satisfied for the presheaf ${{\SD'}}(U)$ the support axiom is automatic for
$\SF_\LOC (U)$. 

To prove the local to global property we first describe the proof for $\SR_\LOC (U)$. Given $A_\a\in\SR_\LOC (U_\a)$ with $A_\a = A_\b$ on $U_{\a\b}$ choose a partition of unity $\{\chi_\a\}$ for $\{ U_\a \}$ with each $\chi_\a$ a characteristic function of a closed set. Then, since $\chi_\a A_\a$ is also a rectifiable current on $U_\a$ (but vanishing near $\p U_\a$), we may consider $\chi_\a A_\a\in \SR_\LOC (X)$ extended by zero to all of $X$.
Set $A= \sum \chi_\b A_\b$ and note that $A\bigm|_{U_\a} = A_\a$. 

The proof of the local to global property for $\SF_\LOC (U)$ can be outlined as  follows.  For a more complete proof see [H,3.1]. Suppose that we are given $T_\a = A_\a + dB_\a$ with $A_\a,B_\a\in \SR_\LOC (U_\a )$
and $T_\a = T_\b$ on $U_{\a\b}$, i.e. $A_\a - A_\b = d(B_\b - B_\a)$. Suppose we could set $T = \sum \chi_\a T_\a$, and verify that 
$$
\chi_\a T_\a = \chi_\a A_\a + \chi_\a dB_\a = \chi_\a A_\a - (d\chi_\a) B_\a + d(\chi_\a B_\a ).\tag $\star$
$$
 More precisely
we must show that $\chi_\a A_\a - d\chi_\a B_\a$ and $\chi_\a B_\a$ define locally rectifiable currents satisfying equation $(\star)$. This is not always true
because  $(d\chi_\a)B_\a$ and/or $\chi_\a dB_\a$ may not be defined. However,
by Federer's theory of slicing, we may choose a slight perturbation of $\chi_\a$ so that $d\chi_\a B_\a$ is a well defined rectifiable current and so that $\chi_\a dB_\a$ is a well defined current with the equation
$d(\chi_\a B_\a ) = (d\chi_\a) B_\a + \chi_\a dB_\a$ satisfied. 

The proof that $\SF_\LOC$ is soft is easier. First consider the analogous result for locally rectifiable currents.  $\SR_\LOC$ is soft because: given a closed set $C\subset X$ and $A\in \SR_\LOC (U)$ where $U$ is a neighbourhood of $C$ we may choose $\chi$ to be the characteristic function of $\overline{V} \subset U$ where $V$ is an open neighbourhood of $C$ and then set $\widehat A = \chi A\in \SR_\LOC (X)$ to be the desired extension. Similarly, given $T = A + dB \in \SF_\LOC (U)$ the current $\widehat T = \chi A  + d(\chi B)$ provides the required extension.
\qed\enddemo

Since $\SF_\LOC$ is soft we have

\proclaim{Corollary A.3} The sheaf $\SF_\LOC$ is acyclic. In particular, for
each locally finite open cover $\SU = \{ U_\a\}$ of $X$, the cohomology $H^1 (\SU, \SF_\LOC) =0$. That is, given $S_{\a\b}\in \SF^p_\LOC (U_{\a\b})$ satisfying
$$ S_{\a\b} \,\, + \,\,S_{\b\g} \,\,+ \,\,S_{\g\a} \,\, = 0  \qquad\text{on } U_{\a\b\g}
$$ 
there exist $S_\a\in\SF^p_\LOC (U_\a)$ such that
$$ S_{\a\b} \,\,=\,\ S_\a \,\,-\,\, S_\b\qquad\text{on } U_{\a\b}.$$

\endproclaim

\proclaim{Theorem A.4}
$$
0\to \bbz \to \SF^0_\LOC\to\SF^1_\LOC\to\dots\to\SF^n_\LOC\to 0$$
is exact.
\endproclaim

\demo{Proof}
Exactness on the left is equivalent to $0\to\bbz \to \SR^0_\LOC\to\SR^1_\LOC$ being exact since $\SF^0_\LOC = \SR^0_\LOC$. This follows immediately from the standard fact that $0\to\bbr \to {{\SD'}}^0\to{{\SD'}}^1$ is exact, i.e. locally a $d$--closed 
generalized function is represented by a constant function. Let $U$ be an open ball in $\bbr^n$. Suppose $T = A + dB$ where
$A\in\SR^p_\LOC (U)$, $B\in\SR^{p-1}_\LOC(U)$ and $p\geq 1$. 
If $T$ is $d$--closed then $dA =0$. Using the standard cone construction (and Federer's theory) there exists $\widehat A\in \SR^{p-1}_\LOC (U)$ such that $d\widehat A = A$ on $U$. Therefore $d(\widehat A + B ) = T$ where $\widehat A + B \in \SR^{p-1}_\LOC (U)\subset \SF^{p-1}_\LOC (U)$. 
\qed\enddemo

\proclaim{Corollary A.5} 
$$
H^p(X,\bbz )\,\, = \,\, \f{\{T\in\SF^p_\LOC (X)  \,\, : \,\, dT =0\}}{d\SF^{p-1}_\LOC (X)}.
$$
 That is, the integrally flat currents  can be used to compute
$\bbz$ cohomology.
\endproclaim

\proclaim{Lemma A.6}
Given $T\in\SF_\LOC (X)$ and a neighbourhood $U$ of $ \oper{spt} T$ there exist $A,B\in \SR_\LOC (X)$ satisfying $T = A + dB$ and with $\oper{spt} A\subset U$, $\oper{spt} B\subset U$. 
\endproclaim

\subheading{Remark}
In particular, Lemma A.6 implies that if $T\in\SF_{\oper{cpt}} (X)$ then there exist  $A,B\in \SR_{\oper{cpt}} (X)$ with $T = A + dB$. Consequently, Federer's definition
of $\SF_{\oper{cpt}} (X)$ (where  $\SF_{\oper{cpt}} (X) = \SR_{\oper{cpt}} (X) + d \SR_{\oper{cpt}} (X)$)
agrees with the definition  given in this appendix (where $\SF_{\oper{cpt}} (X)$ is the space of compactly supported sections of $\SF_\LOC$). 

Federer defines $\SF_\LOC (X)$ to be the space of those currents on $X$ which locally agree with a current of the type $A + dB$ where $A,B \in \SR_{\oper{cpt}} (X)$. As a consequence of the discussion above one can easily show that the definition of $\SF_\LOC (X)$ given in this appendix (namely $\SR_\LOC (X) + d \SR_\LOC (X)$) agrees with the definition in [F].

\demo{Proof of Lemma}
Choose $A,B\in\SR_\LOC (X)$ with $T = A + dB$. On $X\sim \oper{spt} T$, both $B$ and $dB = -A$ are locally rectifiable. By Federer's theory of slicing there exists a neighbourhood $V$ of $\oper{spt} T$ with $\overline V \subset U$ such that on $X\sim \oper{spt} T$, the slice $(d\chi ) B \in \SR _\LOC (X\sim \oper{spt} T )$ exists and $d(\chi B) = (d\chi ) B + \chi dB$ on $X\sim\oper{spt} T$. Here $\chi$ denotes the characteristic function of $\overline V$. Consquently, 
$$ T = A + dB = \chi A + \chi dB = \chi A - (d\chi ) B + d(\chi B),$$
where $\chi A - (d\chi ) B$ and $d(\chi B)$ are locally rectifiable on $X$ with support in $U$.
\qed\enddemo

\subheading{Example 6. Integrally flat $\orX$-twisted  currents}
This example is almost identical to Example 5. Let $\SR^\LOC_k (X)$ denote the space of locally rectifiable  $\orX$-twisted  currents of dimension $k$, and let $\SF_k^\LOC := \SR^\LOC_k (X) + d\SR^\LOC_k (X)$ denote the space of locally integrally flat  $\orX$-twisted  currents of dimension $k$. Then
$$
0\to\wt\bbz_X\to\SF^\LOC_n\to\SF^\LOC_{n-1}\to\dots\to\SF^\LOC_0\to 0
$$
is an acyclic resolution of $\wt\bbz_X$. Hence the complex $\{\SF^\LOC_k,d\} $ can be used to compute $H^p(X,\wt\bbz_X)$, or integer homology $H_k (X,\bbz) = H^p_{\oper{cpt}} (X,\wt\bbz_X)$, where $p+k=n$.
\newline
\newline

\subheading{Example 7} {\bf Integrally flat }  $\SO_V$--{\bf twisted currents.}
Let $V\to X$ be a real bundle with orientation sheaf $\SO_V$ and let $\wt\bbz_V := \bbz\otimes_{\bbz_2}\SO_V$. Let $\wt\SF^p_\LOC (X)$ denote the space of 
degree $p$  
locally integrally flat  $\SO_V$--twisted currents on $X$. Then
$$
0\to \wt\bbz_V\to\wt\SF^0_\LOC\to\wt\SF^1_\LOC\to\dots\to\wt\SF^n_\LOC\to 0
$$ is an acyclic resolution of $\wt\bbz_V$. Hence the complex  $\{ \wt{\SF}^p_\LOC,d\} $  can be used to compute $H^p(X, \wt\bbz_V )$.
\newline
\newline

\subheading{Example 8. Mod q  integrally flat currents}
Cohomology with $\bbz_q$ coefficients can be computed using mod $q$ currents.

\subheading{Definition A.7} The space $\SF^{\oper{mod} q}_\LOC (X)$ of mod $q$ locally integrally flat currents  on $X$  is defined to be the quotient
$\SF_\LOC^p (X) \bigm/ q\SF_\LOC^p (X)$.

To avoid an excess of subscripts and superscripts we always drop the superscript $p$ 
for the degree and the subscript $k$ for the dimension when describing mod $q$ currents.

\subheading{Remark} Federer takes the quotient of $\SF^p_\LOC (X)$ by the closure of $q\SF^p_\LOC (X)$ in the flat topology on $\SF^p_\LOC (X)$ in order to prove compactness theorems for mod $q$ currents. For our purposes --- computing cohomology --- this closure is an unnecessary complication. In addition, Fred Almgren (private communication)  has proven that $q \SF^p_\LOC (X)$ is already closed in $\SF^p_\LOC (X)$. 

\proclaim{Theorem A.8}
The presheaf $\{ \SF^{\oper{mod} q}_\LOC (U) \}$ is a sheaf and this sheaf $\SF^{\oper{mod} q}_\LOC$ is soft.
\endproclaim

\demo{Proof}
Firstly we must verify the support axiom. Suppose $\overline T \in \SF^{\oper{mod} q}_\LOC (X)$ is given (and represented by $T\in \SF_\LOC (X)$). Further suppose $\overline {T}\bigm|_{U_\a} =0$ for each $\a$, i.e. there exist $S_\a\in\SF_\LOC (U_\a )$ such that $T\bigm|_{U_\a} = q S_\a$ for each $\a$. Now $q S_\a = q S_\b$ on $U_{\a\b}$ implies $S_\a = S_\b $ on $U_{\a\b}$. Since $\{ \SF_\LOC (U) \}$ satisfies the local to global axiom there exists $S\in \SF_\LOC (X)$ such that $S\bigm|_{U_\a} = S_\a$. Consequently $T = q S$ on $X$ so that $\overline T = 0$. 

Secondly we must verify the local to global axiom for $\{\SF^{\oper{mod} q}_\LOC (U)\}$.
Suppose $\overline{T}_\a\in \SF^{\oper{mod} q}_\LOC (\Ua )$ are given with $\overline{T}_\a = \overline{T}_\b$ on $U_{\a\b}$. Choosing representatives $T_\a\in \SF_\LOC (U_\a )$ for $\overline{T}_\a$, the equation $\overline{T}_\a = \overline{T}_\b$ on $U_{\a\b}$ says that there exist $S_{\a\b}\in \SF_\LOC (U_{\a\b})$ such that $T_\a - T_\b = q S_{\a\b}$ on $U_{\a\b}$. The cocycle condition for $S_{\a\b}$ is satisfied since it is satisfied for $q S_{\a\b}$. Recall that $H^1 (\SU, \SF_\LOC ) =0$. This implies that there exist $S_\a\in\SF_\LOC (\Ua )$ such that $S_\a - S_\b = S_{\a\b}$. Consequently 
$$T = T_\a - q S_\a\in\SF_\LOC (X)$$
 is well defined independent of $\a$ and provides the global representation of a class $\overline T\in \SF^{\oper{mod} q}_\LOC (X)$ with $\overline{T}\bigm|_{U_\a} = \overline{T}_\a$ for each $\a$.

Finally we must show that the sheaf $\SF^{\oper{mod} q}_\LOC$ is soft. Suppose $\overline T \in \SF^{\oper{mod} q}_\LOC (U)$ is given (represented by $T\in \SF_\LOC (U)$) and $U$ is an open neighbourhood of a closed set $C$. Since $\SF_\LOC$ is soft there exists $S\in \SF_\LOC (X)$ such that $S\bigm|_V = T\bigm|_V$ for some neighbourhood $V\subset U$ of $C$. Now $\overline S \in \SF^{\oper{mod} q}_\LOC (X)$ and $\overline S = \overline T$ on the neighbourhood $V$ of $C$.
\qed\enddemo

\proclaim{Proposition A.9}
Given $\overline T \in \SF^{\oper{mod} q}_\LOC (X)$ and a neighbourhood $U$ of $A := \oper{spt} \overline T$ there exists a representative $T\in\SF_\LOC (X)$ of $\overline T$ with $\oper{spt} T \subset U$.
\endproclaim

\demo{Proof}
Let $T\in\SF_\LOC (X)$ denote an arbitrary representative of $\overline T$. Restricted to $X\sim A$, $T\bigm|_{X\sim A} = q \widehat S$ for some $\widehat S \in \SF_\LOC (X\sim A )$. Since $\SF_\LOC$ is soft there exists a global section $S\in \SF_\LOC (X)$ which agrees with $\widehat S$ on a neighbourhood of the closed set $X\sim U$. Now 
 $T - q S\in\SF_\LOC (X)$ is another
representative of $\overline T$  and  $T - q S$ vanishes on a neighbourhood of $X\sim U$.\qed\enddemo

Note that the sheaf  $\bbz_q$  is a subsheaf of $\SF^{\oper{mod} q}_\LOC$ when the degree is zero (i.e. when dimension is $n$).

\proclaim{Lemma A.10}
Suppose $\overline T\in \SF^{\oper{mod} q}_\LOC (U)$ is of degree zero on the open unit ball $U$. If $d\overline T = 0$ then $\overline T$ is represented by a constant integer valued function.
\endproclaim

\demo{Proof}
The equation $d\overline T = 0$ means that $dT = q S$ with $S\in \SF^1_\LOC (U)$  and $T\in \SF^0_\LOC (U)$  a representative for $\overline T$. The current $q S$ and hence $S$ is $d$--closed. Therefore there exists $R\in\SF^0_\LOC (U)$  satisfying $dR = S$ (because this is true for locally rectifiable currents). Now $T- q R\in \SF^0_\LOC (U)$ also represents $\overline T$ and is $d$--closed. Therefore $T - q R$ is a constant integer
valued function.
\qed\enddemo

This proves that $0\to {\bbz_q}\to \SF^{\oper{mod} q}_{\LOC, \oper{deg} 0}\overset d \to \to \SF^{\oper{mod} q}_{\LOC, \oper{deg} 1}$ is exact.

\proclaim{Lemma A.11}
Let $p\geq 1$. Given $\overline T \in   \SF^{\oper{mod} q}_{\LOC, \oper{deg} p} (U)$
on the unit ball $U$ satisfying $d\overline T = 0$ there exists $\overline S \in  \SF^{\oper{mod} q}_{\LOC, \oper{deg} p-1} (U)$ with $d\,\overline S = \overline T$.
\endproclaim

\demo{Proof}
There exist $T\in \SF^p_\LOC (U)$ and $S\in\SF^{p+1}_\LOC (U)$ such that $dT = q S$. Thus $dS =0$. Now $S = A + dB$ where $A,B\in\SR_\LOC (U)$. Thus 
$dA =0$. Since degree $A = p+1 \geq 1$ there exists $R\in\SR^p_\LOC (U)$ such that $A = dR$. Therefore $S = dB$ for some $B\in \SR^p_\LOC (U)$. Consequently
$\widehat T = T - q B\in \SF^p_\LOC (U)$ is another representative of $\overline T$ 
with $d\widehat T = 0$. Finally solve $dS = \widehat T$ on $U$ with $S\in \SF^p_\LOC (U)$.
\qed\enddemo

In summary, the mod $q$ currents  $\SF^{\oper{mod} q}_\LOC (X)$  may be used to compute cohomology with ${\bbz_q}$ coefficients.

\proclaim{Theorem A.12}$$
H^p (X,{\bbz_q} ) \,\, = \,\, \f{ \{\overline T \in  \SF^{\oper{mod} q}_{\LOC, \oper{deg} p} (X)\,\, : \,\, d\overline T = 0\} }{ d\,  \SF^{\oper{mod} q}_{\LOC, \oper{deg} p-1} (X) }.
$$
\endproclaim

\subheading{Example 9}{\bf  Mod q integrally flat } $\orX${\bf -twisted
 currents.}
This example is almost identical to Example 8. Let $\SF^\LOC_{\oper{mod} q} (X)
:= \SF^\LOC_k (X) \bigm/ q \SF^\LOC_k (X)$ be the space of mod $q$ locally integrally flat $\orX$-twisted  currents on $X$ and let $\wt\bbz_q := \SO_X\otimes_{\bbz_2}\bbz_q$. Then the complex  $\{ \SF^\LOC_{\oper{mod} q} (X)\, , \, d\}$ can be used to compute $H^p(X,\wt\bbz_q)$, or $H_k(X,\bbz_q) = H^p_{\oper{cpt}} (X,\wt\bbz_q)$ where $p+k=n$.

\subheading{Example 10. The mod 2 integrally flat case}
In this example we wish to compute $\bbz_2$ cohomology $H^p(X,\bbz_2)$. This can be done in several equivalent ways. Firstly, restating Theorem A.12 in the case $q=2$, $H^p(X,\bbz_2)$ can be computed from the complex $\SF^{\oper{mod} 2}_\LOC  (X) := \SF^*_\LOC (X) \bigm/ 2 \SF^*_\LOC (X)$ of mod 2 locally integrally flat currents  on $X$.

Next suppose that $V$ is a real bundle with orientation sheaf $\SO_V$ and let 
$\wt\SF^p_\LOC (X)$ denote the  degree $p$ locally integrally flat $\SO_V$--twisted currents  on $X$. 

\proclaim{Lemma A.13}
$$
 \SF^p_\LOC (X)\bigm/ 2\SF^p_\LOC (X)\,\, \cong \,\,\wt\SF^p_\LOC (X)\bigm/ 2\wt\SF^p_\LOC (X)
$$
Therefore $H^p(X,\bbz_2)$ can also be computed using the complex of mod 2
locally integrally flat  $\SO_V$--twisted currents.
\endproclaim

\demo{Proof of Lemma}
First we define a mapping $\varphi : {\wt\SF}^p_\LOC (X) \to \SF^p_\LOC (X) \bigm/ 2 \SF^p_\LOC (X)$. Let $T\in \wt{\SF}^p_\LOC (X)$. Choose a locally finite open cover $\{ U_\a \} $ of $X$ and fix  ordered frames $e_\a$ for $V$ over $U_\a$. For each pair
$(U_\a, e_\a)$ let $T_\a\in \SF^p_\LOC (U_\a )$ be defined by $T = [e_\a]\otimes T_\a$ on $U_\a$. Note that $T_\a = \pm T_\b$ on $U_a\cap U_b$.
 Define $S_{\a\b}\in\SF_{\oper{loc}} (U_\a\cap U_\b )$ by
$$
S_{\a\b} = \cases 0 &\text{if  $T_\a = T_\b$ on $U_\a\cap U_\b$,}\\
                  T_\a &\text{if $T_\a = \, -T_\b$  on $U_\a\cap U_\b$.}\endcases
$$
Then 
$$
T_\a \,\, - \,\, T_\b \,\, = \,\, 2 S_{\a\b} \qquad\text{on } U_\a\cap U_\b.
\tag A.14
$$
Let $\overline{T_\a}\in \SF^{\oper{mod} 2}_\LOC (U_\a )$ be the mod 2 class of $T_\a$. Then, by equation (A.14), $\overline{T_\a} = \overline{T_\b}$ on $U_\a\cap U_\b$. Since the presheaf $\{ \SF^{\oper{mod} 2}_\LOC (U )\} $ satisfies the Local to Global axiom (see Theorem A.8)
there is a well defined element $\varphi (T) \in \SF^{\oper{mod} 2}_\LOC (X )$
so that $\varphi (T)\bigm|_{U_\a} = \overline{T_\a}$. Note that $\varphi (T)$ is well defined independent of the choices of locally finite open cover $\{ U_\a \} $ and frames $e_\a$ for $V$ over $U_\a$. 

In summary we have defined a map 
 $\varphi : {\wt\SF}^p_\LOC (X) \to \SF^p_\LOC (X) \bigm/ 2 \SF^p_\LOC (X)$.
Now, since $\varphi (2T)\bigm|_{U_\a} = \overline{2 T_\a} = 0$, the induced map $\widehat\varphi : {\wt\SF}^p_\LOC (X)\bigm/ 2  {\wt\SF}^p_\LOC (X) \to \SF^p_\LOC (X) \bigm/ 2 \SF^p_\LOC (X)$ is well defined. Finally, to show that $\widehat\varphi$ is an isomorphism, we can use the same procedure to construct an inverse for $\widehat\varphi$.
\qed\enddemo

\enddocument